\documentclass[12pt]{article}
\usepackage{amsmath}
\usepackage{amssymb}
\usepackage{amsthm}
\usepackage{amscd}
\usepackage{amsfonts}
\usepackage{graphicx}
\usepackage{fancyhdr}
\usepackage{color}
\usepackage{enumerate}
\usepackage{amsfonts}
\usepackage{amsthm}
\usepackage{psfrag} 
\usepackage{epsfig}
\usepackage{subfigure}
\usepackage{graphicx}
\usepackage{amssymb}
\usepackage{epstopdf}
\usepackage{amsbsy}
\usepackage{latexsym}
\usepackage{moreverb}                      
\usepackage{textcomp}

\usepackage{multicol}

\usepackage{algorithm}
\usepackage{algorithmic}

\usepackage[top=2.8cm,bottom=2.8cm,left=2.5cm,right=2.5cm]{geometry}
\usepackage[T1]{fontenc}
\usepackage[utf8]{inputenc}

\usepackage[affil-it]{authblk} 
\usepackage{etoolbox}
\usepackage{lmodern}

\makeatletter
\patchcmd{\@maketitle}{\LARGE \@title}{\fontsize{18}{19.2}\selectfont\@title}{}{}
\makeatother

\usepackage{sectsty}
\sectionfont{\fontsize{13}{15}\selectfont}
\subsectionfont{\fontsize{12}{15}\selectfont}

\newcommand{\daniel}[1]{{{\color{black} {#1}}}}

\begin{document}

\title{Wasserstein metric-driven Bayesian inversion\\
with applications to signal processing}

\author[1]{M. Motamed\thanks{motamed@math.unm.edu}}
\author[2]{D. Appel{\"o}\thanks{daniel.appelo@colorado.edu}}
\affil[1]{Department of Mathematics and Statistics, University of New Mexico, USA}
\affil[2]{Department of Applied Mathematics, University of Colorado Boulder, USA}

\date{\small \today}

\maketitle

\medskip

{\small
\noindent
{\bf Abstract.} 
We present a Bayesian framework based on a new exponential
likelihood function driven by the quadratic Wasserstein metric. Compared
to conventional Bayesian models based on Gaussian likelihood functions
driven by the least-squares norm ($L_2$ norm), the new framework features several
advantages. First, the new framework does not rely on the likelihood
of the measurement noise and hence can treat complicated noise
structures such as combined additive and multiplicative noise. Secondly, unlike the normal likelihood function,
the Wasserstein-based exponential likelihood function does not usually
generate multiple local extrema. As a result, the new framework
features better convergence to correct posteriors when a
Markov Chain Monte Carlo sampling algorithm is employed. Thirdly, in the
particular case of signal processing problems, while a normal likelihood function
measures only the amplitude differences between the observed and
simulated signals, the new likelihood function can capture both 
amplitude and  phase differences. We apply the new framework to a
class of signal processing problems, that is, the
inverse uncertainty quantification of waveforms, and demonstrate its advantages compared to Bayesian models with normal
likelihood functions.

}


\medskip
\medskip


\section{Introduction}

The Euclidean distance 
is one of the most widely used metrics in constructing misfit functions in
deterministic inversion algorithms. In a Bayesian framework (see e.g. \cite{Gelman_etal:04,Kaipo_Somersalo:05,Stuart:10}), this misfit
function corresponds to a normal likelihood function, which is valid
under the (often unrealistic) assumption that the experimental noise
is additive and Gaussian. In addition to being unrealistic, a Gaussian
likelihood function measures only the amplitude differences between the
observed and simulated signals and does not capture the phase
differences. Moreover, Gaussian likelihood functions often have
multiple local extrema that may result in converging to wrong
posteriors (see Section \ref{sec:num1}) when a Markov Chain Monte
Carlo (MCMC) sampling algorithm \cite{MCMC:04} is used. 
Such disadvantages are particularly of great concern in signal
processing applications, such as seismic inversion problems, where the
noise is not necessarily additive and the wave signals are oscillatory and cyclic.

In the present work we will propose a Bayesian framework based on an exponential
likelihood function driven by the quadratic Wasserstein metric; see
e.g. \cite{Villani:03,Villani:09}. 
We will compare the new framework with current Bayesian models based
on normal likelihood functions when applied to a class of signal
processing problems, including seismic source inversion and seismic
imaging. 
\daniel{In particular the ideas presented here find their inspiration in the recent activity in the area of (deterministic) Full Waveform Inversion (FWI) exploiting distances connected to optimal transportation. The use of the quadratic Wasserstein distance in FWI is relatively new but has had substantial impact in the exploration community as well as in the mathematical geosciences community. Engquist and Froese first proposed to use the quadratic Wasserstein distance for seismic signals in \cite{engquist2014application}.  This work has been analyzed and extended in a sequence of papers by Yang, Engquist,  and co-authors \cite{engquist2016optimal,Wasserstein_Bjorn,2018arXiv180804801E} clearly demonstrating the power of the approach when applied to FWI. We note that the quadratic Wasserstein metric has recently been used to invert for source mechanisms and location in \cite{ChenCWY18,2018arXiv181001710B}. 

We further note that in FWI there is also considerable interest in the approach 
by M{\'e}tivier and co-authors, \cite{metivier2016measuring,metivier2016optimal,metivier2018graph,FrenchOT3,FrenchOT4}. In their approach the quadratic Wasserstein distance is replaced by a relaxed 1-Wasserstein distance obtained from the Monge-Kantorovich problem via a maximization in the space of bounded 1-Lipschitz functions. Unlike the quadratic Wasserstein distance the resulting distance does not rely on normalizing the input signals to be positive and with equal mass. 

Given the success of both approaches in FWI, a natural extension of this work would be to use the techniques from \cite{metivier2016measuring,metivier2016optimal,metivier2018graph,FrenchOT3} in place of the quadratic Wasserstein distance. 

In this paper we focus on the quadratic Wasserstein and extend the work \cite{engquist2014application,engquist2016optimal,Wasserstein_Bjorn,2018arXiv180804801E} to Bayesian inversion and demonstrate applications of the new framework in the uncertainty quantification of seismic inversion problems.
}


The rest of the paper is organized as follows. In Section \ref{sec:bayes} we
formulate the Bayesian inversion problem and discuss the likelihood
and measurement noise structures. We then present the new
Wasserstein metric-driven Bayesian framework in Section \ref{sec:W2bayes}, followed by a standard
MCMC algorithm in Section \ref{sec:MCMC}. We finally present several numerical
examples related to seismic wave inversion in Section \ref{sec:num1}.

\section{Bayesian inversion problem} \label{sec:bayes}

\subsection{Problem formulation}

Let ${\bf g} = (g_1, \dotsc, g_N) \in {\mathbb R}^{N}$ be a vector of $N$ observed
quantities, for instance, measured by an experimental device. Let
further ${\bf f} = (f_1, \dotsc, f_N) \in {\mathbb R}^{N}$ be a vector of $N$ predicted
quantities corresponding to ${\bf g}$, computed by a forward
predictive model with an input parameter vector $\boldsymbol\theta =
(\theta_1, \dotsc, \theta_m) \in \Theta  \subset {\mathbb R}^m$:
$$
{\bf f} = {\bf f}(\boldsymbol\theta): \Theta  \subset {\mathbb R}^m \rightarrow {\mathbb R}^{N}.
$$
A main assumption in mathematical statistics is that the model
parameter vector $\boldsymbol\theta$, and consequently any physical quantity
to be predicted by the forward model, is described by a stochastic law
characterized by a probability density function (PDF). The density of
$\boldsymbol\theta$ is available to us only through the experimental measurements ${\bf g}$ that are often corrupted
by experimental noise and error. The goal is then to find the density of $\boldsymbol\theta$, given the forward model ${\bf f}(\boldsymbol\theta)$ and the noisy data ${\bf g}$. 

In Bayesian inversion the conditional posterior probability density of
the model parameter vector $\boldsymbol\theta$ is given by the Bayes'
rule; see e.g. \cite{Bayes:1763,Gelman_etal:04,Kaipo_Somersalo:05,Stuart:10}:
$$
\pi(\boldsymbol\theta | {\bf g}) = \frac{\pi({\bf
    g} | \boldsymbol\theta) \, \pi(\boldsymbol\theta)}{ \int_{\Theta}  \pi({\bf
    g} | \boldsymbol\theta) \, \pi(\boldsymbol\theta) \, d \boldsymbol\theta}
\propto \pi({\bf
    g} | \boldsymbol\theta) \, \pi(\boldsymbol\theta),
$$
where $\pi({\bf g} | \boldsymbol\theta)$ is the
likelihood, and $\pi(\boldsymbol\theta)$ is the prior density of
$\boldsymbol\theta$. The proportionality follows from the fact that the
denominator (referred to as the evidence) does not depend on
$\boldsymbol\theta$ and hence can be regarded as a constant. 
Assuming the observed samples $\{ g_i \}_{i=1}^N$ are
independent, their joint likelihood reads
$$
 L(\boldsymbol\theta) := \pi({\bf g} | \boldsymbol\theta) = \prod_{i=1}^N \pi(g_i |
\boldsymbol\theta),
$$
where the first equality indicates that the likelihood can be viewed as a
function of only $\boldsymbol\theta$, referred to as the {\it likelihood
function} and denoted by $L(\boldsymbol\theta)$.

\subsection{Likelihood and noise structure}

A major step in Bayesian analysis is the selection of the
likelihood function. This is usually done based on the structure of
the measurement noise. 
The most common choice is the Gaussian likelihood,
built based on the assumption that the measurement noise $\{
\varepsilon_i \}_{i=1}^N$ in the $N$
measured quantities $\{ g_i \}_{i=1}^N$ is additive and normally
distributed with zero mean and a standard deviation
$\sigma$: 
$$ 
g_i = f_i(\boldsymbol\theta) + \varepsilon_i, \qquad \varepsilon_i \sim
\text{Normal}(0, \sigma), \qquad i=1, \dotsc, N.
$$ 
In this case the likelihood function reads
$$
L_{\text{norm}}(\boldsymbol\theta) =\pi_{\text{norm}}( {\bf g} | \boldsymbol\theta) = \frac{1}{(2 \, \pi)^{N/2} \,
  \sigma^N} \exp \bigl( \frac{-1}{2 \, \sigma^2}  \sum_{i=1}^{N} |
 g_i - f_i(\boldsymbol\theta)|^2    \bigr).
$$

One critical problem with this setting is that the Gaussian additive
noise assumption may not be realistic. For
instance, in ultrasound and laser imaging with ${\bf g} \in {\mathbb R}^N$ being
an $N$-pixel image, the noise may be multiplicative. A practical
example is a speckle noise (see e.g. \cite{Speckle:07}) with mean one and variance $1/s$:
$$ 
g_i = \varepsilon_i \, f_i(\boldsymbol\theta), \qquad \varepsilon_i \sim
\text{Gamma}(s, 1/s), \qquad i=1, \dotsc, N.
$$ 
Here $\text{Gamma}(a,b)$ is a Gamma distribution with a shape parameter $a>0$
and a rate parameter $b>0$. In some applications, the noise structure
may be even more complicated. For instance the measurement noise may
be both additive and multiplicative, taking the form
\begin{equation}\label{noise_structure}
g_i = \varepsilon_i^{(1)} \, f_i(\boldsymbol\theta) +
\varepsilon_i^{(2)} , \qquad i=1, \dotsc, N.
\end{equation}
Such complicated noise structures may not be easily treated by Bayesian frameworks that rely on the likelihood of the
noise, even when a non-Gaussian likelihood function is
employed. It is to be noted that the noise structure is not the only
problem with Gaussian likelihoods. 
Even if the noise structure is additive and can be modeled by a Gaussian likelihood function, the
likelihood may have multiple local extrema, resulting in wrong
posteriors when applying a MCMC sampling method; see Section
\ref{sec:num1}. Moreover, in signal processing problems, a Gaussian
likelihood function that is based on the $L_2$ norm measures only the
amplitude differences between the observed and predicted signals and
does not capture the phase differences. This latter problem is indeed a crucial
limitation for instance in: a) seismic inversion with moderately to highly oscillatory seismic waves, and b) MRI noise reduction where the retainment of the image fine features is desired.

In what follows we present a Bayesian framework based on an exponential
likelihood function driven by a quadratic Wasserstein metric. Unlike conventional
Bayesian analysis, this framework does not rely on the likelihood of the measurement
noise and hence can treat complicated noise structures as in \eqref{noise_structure}. Moreover,
since the Wasserstein metric has better optimization and fitting
properties than the least-squares norm (see
  e.g. \cite{Wasserstein_Bjorn} and the references therein), the new framework may feature
better convergence properties than the standard frameworks based on Gaussian likelihoods.

\section{Wasserstein metric-driven Bayesian framework} \label{sec:W2bayes}

The Wasserstein metric is a distance function defined between two
probability distributions. It corresponds to the minimum
``cost'' of turning one distribution into the other. The metric has been applied to various fields, including image processing \cite{PapHab},
computer vision \cite{7780937}, Stochastic programming, \cite{MohajerinEsfahani2017}, and seismic imaging \cite{Wasserstein_Bjorn}. We
refer to \cite{Villani:03,Villani:09} for a complete introduction to
Wasserstein distances. 

Consider two discrete-time signals ${\bf f}, {\bf g} \in {\mathbb
  R}^N$ given at a discrete set of time level $\{ t_i
\}_{i=1}^{N}$. 
{\color{black}
Since the two (oscillatory) signals are not necessarily probability
distributions, i.e. positivity and mass balance are not expected, the
signals need to be pre-processed before we can measure their Wasserstein
distance. Here, 
we follow the normalization via a linear transformation and rescaling
proposed in \cite{Wasserstein_Bjorn} and construct two discrete CDFs as follows:
}

\begin{itemize}
\item Shift the signals to ensure positivity: select a constant $c$ so
  that $f_i+ c > 0$ and $g_i+c >0$, for all $i = 1, \dotsc, N$.

\item Rescale the shifted signals so that they share a common unit total mass:
\begin{equation}\label{tilde_fg}
\tilde{\bf f} = \frac{{\bf f}+ c }{ \langle {\bf f} + c \rangle},
\qquad \tilde{\bf g} = \frac{{\bf g}+ c }{ \langle {\bf g}+ c \rangle},
\end{equation}
where $\langle {\bf f} \rangle = \sum_{i=1}^N f_i$ is the total mass of the signal ${\bf f}$.

\item Find the discrete CDFs $F=(F_1, \dotsc, F_N)$ and $G=(G_1, \dotsc, G_N)$ with the components 
\begin{equation}\label{FG}
F_i = \sum_{j=1}^i \tilde{f}_j, \qquad G_i = \sum_{j=1}^i \tilde{g}_j,
\qquad i=1, \dotsc, N,
\end{equation}
where $\tilde{f}_j$ and $\tilde{g}_j$ are the $j$-th components of
$\tilde{\bf f}$ and $\tilde{\bf g}$ in \eqref{tilde_fg}, respectively.

\end{itemize}

The discrete quadratic Wasserstein distance between the two signals then reads
\begin{equation} \label{W2_distance}
d_W({\bf f},{\bf g}) = \sum_{i=1}^{N} |  t_i -
T_i |^2 \,  \tilde{f}_i,
\qquad T=G^{-1} \circ F.
\end{equation}
We note that $T=G^{-1} \circ F$ is the ``optimal'' map from the density
$\tilde{\bf f}$ to the density $\tilde{\bf g}$. Intuitively, if we view the
two densities as two piles of dirt with the same unit mass, the map
would correspond to the minimum ``cost'' of turning one pile into the
other, i.e. the minimum amount of dirt that needs to be moved times
the distance it has to be moved. 
This corresponds to the horizontal distance between the two discrete CDFs $F$
and $G$. For computing the metric $d_W$ in \eqref{W2_distance}, the map $T$ is computed by interpolation.

We then consider an exponential likelihood function 
\begin{equation} \label{W2_likelihood}
L_{\text{exp}}(\boldsymbol\theta) = \pi_{\text{exp}}( {\bf g} | \boldsymbol\theta) =s^N \, \exp ( -s
\, d_W({\bf f}(\boldsymbol\theta),{\bf g})),
\end{equation}
where $s>0$ is a rate parameter and will be considered as a
hyperparmeter in the MCMC sampling algorithm. {\color{black}The main motivation for
the selection of an exponential function is that we would like to
recover the Wasserstein-driven cost functional in deterministic
inverse problems by taking the negative logarithm of the posterior in the
Bayesian setting \cite{Stuart:10}.}

It is to be noted that the Wasserstein metric can be computed in any
dimension through the (numerical) solution  of an Monge--Amp\`{e}re
equation, see e.g. \cite{BENAMOU2014107}. For dimensions larger than
one, where there are more than one signal, this may be a
costly procedure. We therefore follow the work by Yang
et. al. \cite{Wasserstein_Bjorn} and use a ``trace-by-trace''
approach, summing up one dimensional quadratic Wasserstein distances;
see Section \ref{sec:num1} where there are more than one signal
collected by several receivers.

\begin{figure}[]
\begin{center}
\includegraphics[width=0.32\textwidth]{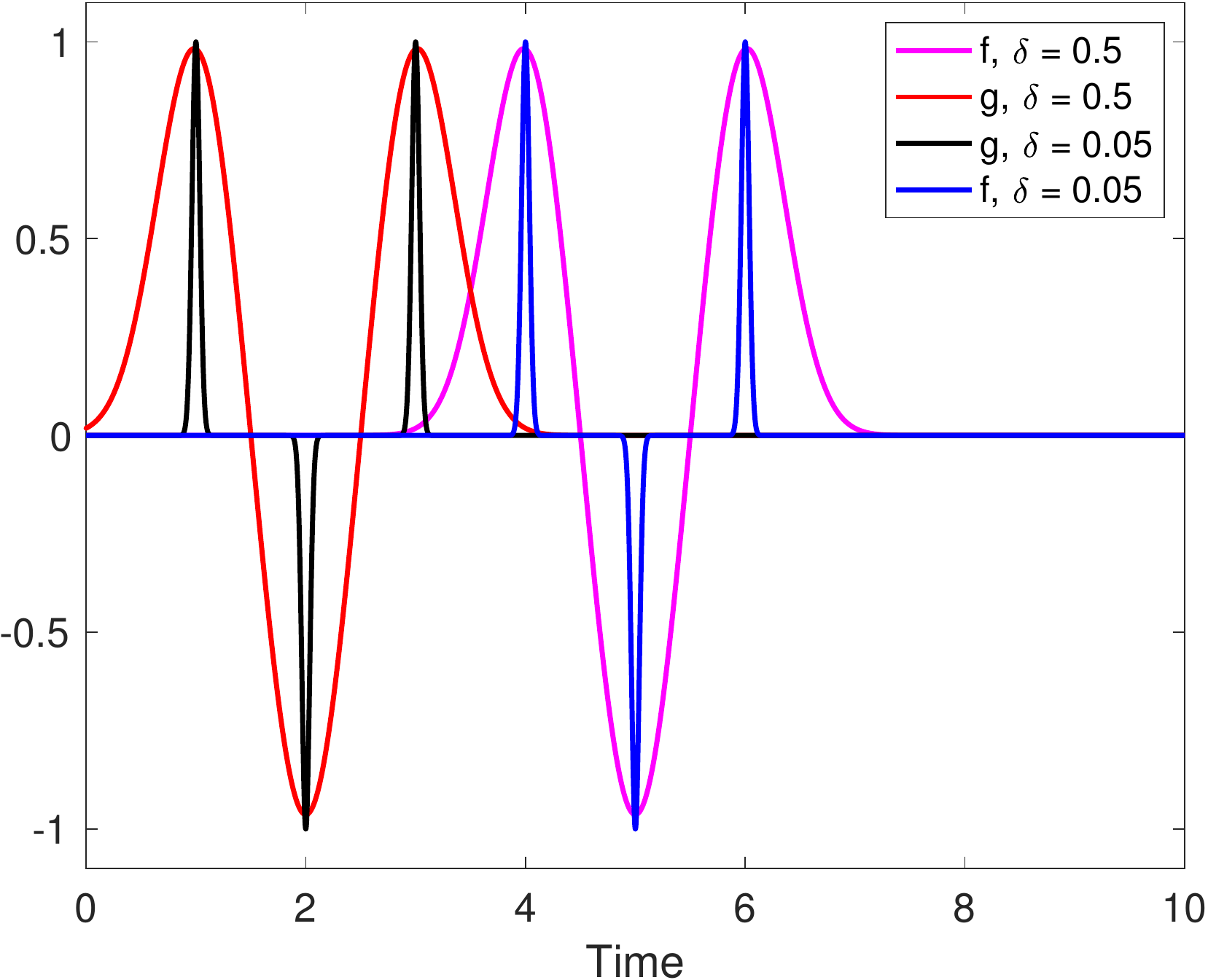}  
\includegraphics[width=0.32\textwidth]{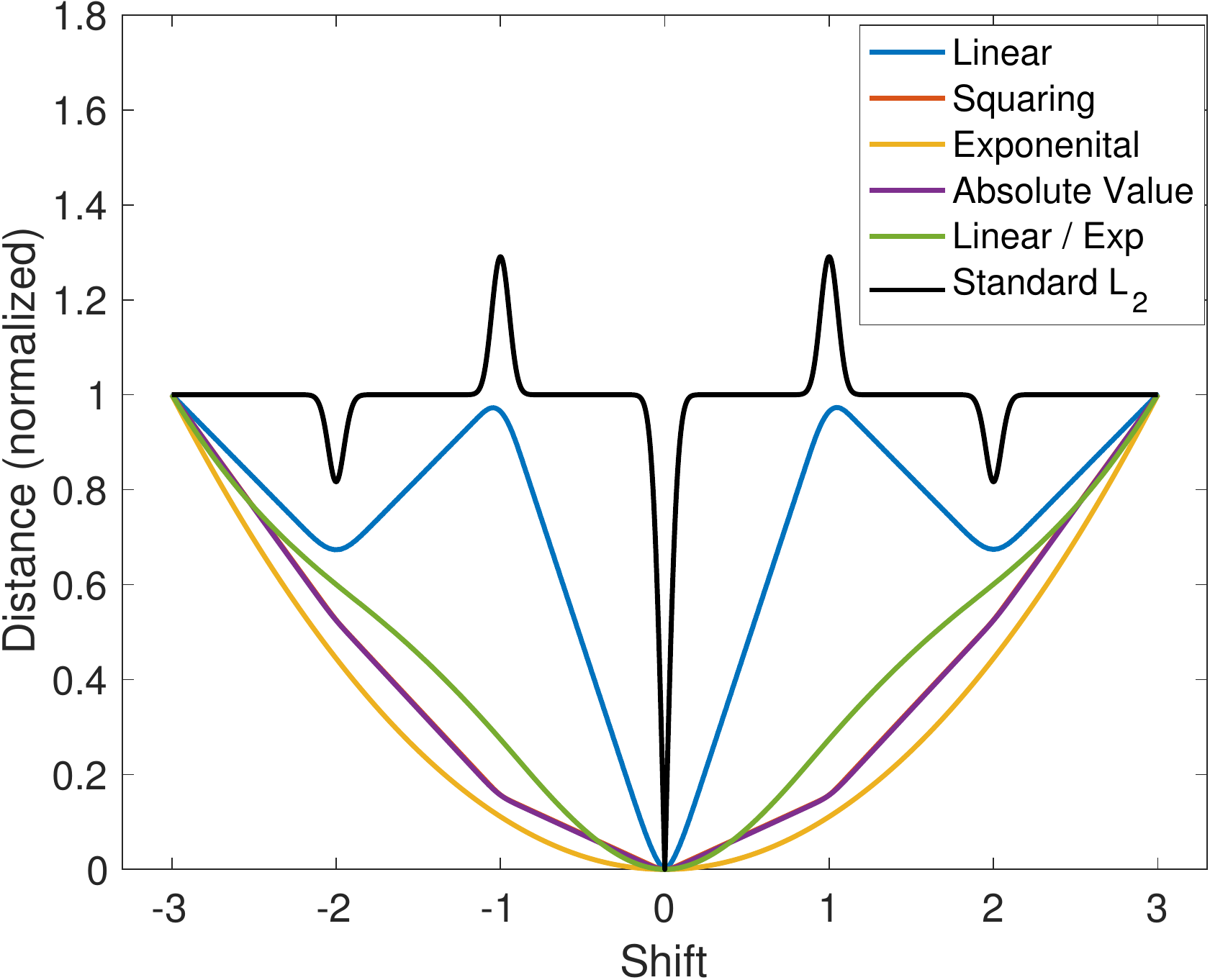}      
\includegraphics[width=0.32\textwidth]{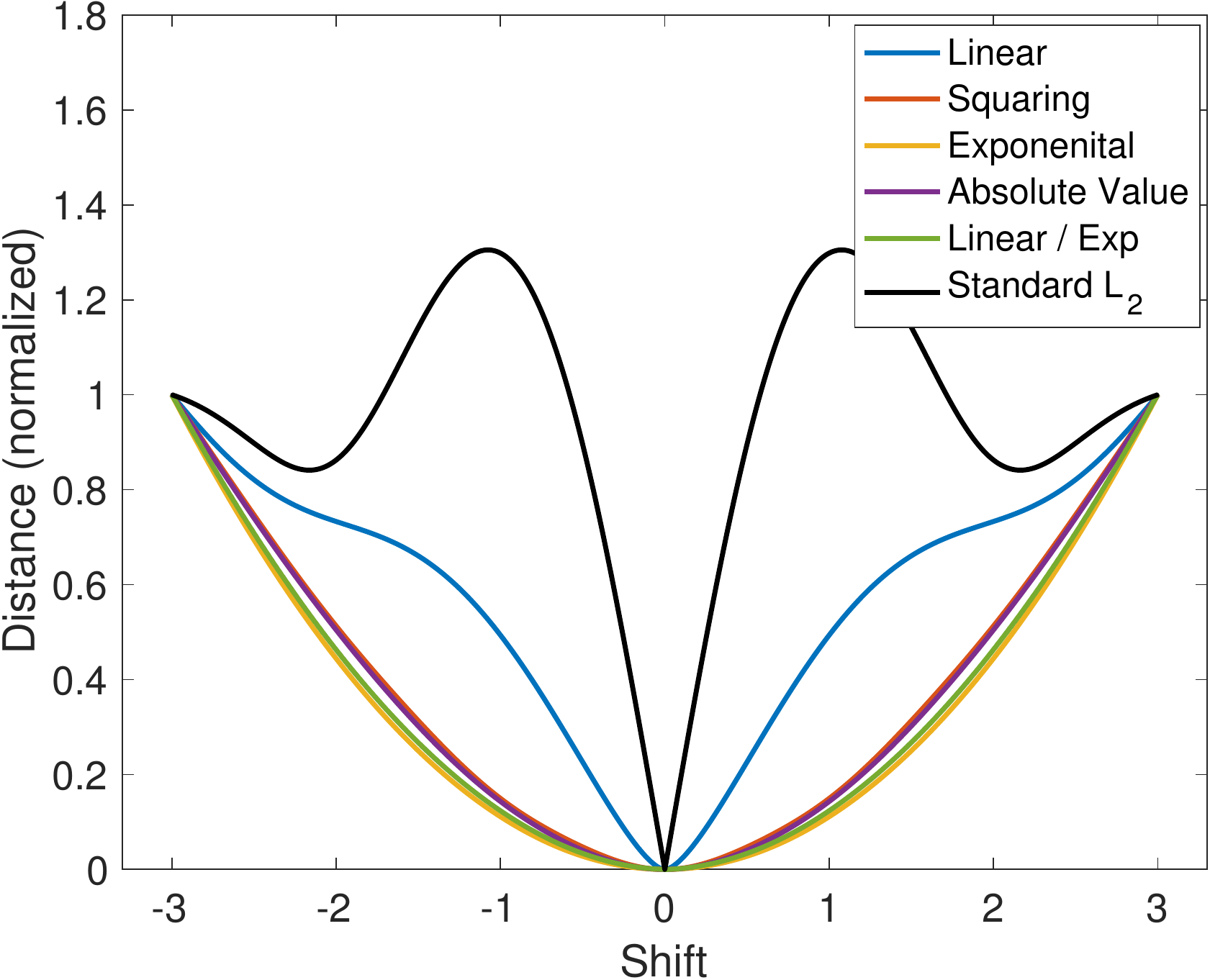}     
\caption{Comparison of different normalization / scaling procedures for the quadratic Wasserstein distance and the standard $L_2$ distance for two signals with different width (see the text for descriptions of the legend). The middle figure is for the narrow signals and the right is for the wide signals. \label{fig:W2_norm_comp}}
\end{center}
\end{figure}

\medskip
\noindent
{\bf Convexity of the Wasserstein exponential likelihood.} 
{\color{black}
One important feature of the quadratic Wasserstein metric is its
convexity with respect to the phase shift, phase dilation, and amplitude change in the simulated and measured signals
\cite{engquist2016optimal}. 
However, due to normalization during the pre-processing step, the
convexity may be lost. This may in turn lead to the cycle skipping
issue in the proposed Bayesian framework, i.e. the presence of
multiple local optima in the exponential likelihood function
\eqref{W2_likelihood}. 
We notice that various normalization strategies, in the context of FWI, have been proposed to
overcome this issue; see e.g. \cite{engquist2014application,engquist2016optimal,Wasserstein_Bjorn}. 
To see the effect of different normalization strategies on the
  convexity of the quadratic Wasserstein distance with respect to
  phase shifts and to compare the cycle skipping issue with that of
  the usual $L_2$ distance, we consider a simple example with two identical but shifted signals. Precisely the original signal $f$ is  
\[
f(t) =  e^{-\left( \frac{t-4}{\delta}\right)^2} 
- e^{-\left( \frac{t-5}{\delta}\right)^2} 
+e^{-\left( \frac{t-6}{\delta}\right)^2},
\] 
and $g$ is a shifted version of $f$
\[
g(t) =  e^{-\left( \frac{t-s-4}{\delta}\right)^2} 
- e^{-\left( \frac{t-s-5}{\delta}\right)^2} 
+e^{-\left( \frac{t-s-6}{\delta}\right)^2},
\] 
with shift $s$. In Figure \ref{fig:W2_norm_comp} we compare the convexity of the distance between $f$ and $g$ as a function of the shift and for different normalizations. The values of the distances have been scaled to start at one when the shift is -3.} 
\daniel{The different options to preserve positivity we use are }
\begin{eqnarray*}
&\hat{f} = f + c, & \text{ Linear scaling,}\\
&\hat{f} = f^2, & \text{ Square scaling,}\\
&\hat{f} = \exp(cf), & \text{ Exponential scaling,}\\
&\hat{f} = |f|, & \text{ Absolute value scaling,}\\
&\hat{f} = \left\{ \begin{array}{c}
\exp(cf) \text{ if } f<0, \\
f+1/c \text{ if } f \geq 0, \\
\end{array} \right. & \text{ Linear/Exponential scaling,}
\end{eqnarray*}  
\daniel{The generic constant $c$ should be chosen big enough to ensure positivity for the linear scaling and small enough for exponential scaling not to overflow. 

The results in Figure \ref{fig:W2_norm_comp} show that, in general, one should not expect to preserve the convexity of the quadratic Wasserstein distance when the linear scaling is used, as has been discussed recently in e.g. \cite{2018arXiv180804801E} and \cite{FrenchOT3}. However, if the signals are not clearly separated (corresponding to a larger $\delta$) the linear scaling does preserve the convexity, see the rightmost subfigure. Compared to the $L_2$ distance the basin of attraction is significantly expanded even in the case with local minima. The other alternatives appears to give quite good results in terms of convexity, but it should be pointed out that the squaring and absolute value functions are not unique, in the sense that many signals $f$ can produce $\hat{f}$ and that these options are therefore typically not used in deterministic problems. The most robust option is probably the combination of the exponential and linear procedure.     

We note that in all our computations below we use the linear scaling. For the one dimensional examples the main signals are positive but the added noise can make them slightly negative, activating the normalization procedure. However as the size of the noise is small we expect that the linear scaling results in a convex cost function. In the two dimensional examples below the discrete signals typically contain equal parts of positive and negative entries so that we cannot exclude the possibility of local minima. We plan to carry out a detailed exploration of the impact of the scaling on the proposed method in the future.}  

\section{Numerical algorithm} \label{sec:MCMC}

We employ a Metropolis-Hastings-within-Gibbs sampling strategy (see
e.g. \cite{Gelman_etal:04}) within the proposed Bayesian
framework. The algorithm, summarized in Algorithm \ref{ALG_MG},
iteratively generates a sequence of samples, forming a Markov chain,
whose distribution approaches the target distribution of parameters in
the limit. A candidate sample is produced based on the current sample
value and is later accepted or rejected with some probability $\alpha$; see below. 
The algorithm consists of two interactive parts: a Gibbs sampler and a Metropolis-Hastings sampler. The Gibbs sampler is used to sample the
hyperparameter $s>0$ with a fixed parameter vector
$\boldsymbol\theta$, and the Metropolis-Hastings sampler is used to
sample $\boldsymbol\theta$ with a fixed hyperparameter $s$.

\medskip
\noindent
{\bf Gibbs sampler.} 
The Gibbs sampler employs the proposed exponential likelihood
\eqref{W2_likelihood} and a Gamma prior on its
rate parameter $s$, i.e. $s \sim \text{Gamma}(a,b)$, 
where $a>0$ and $b>0$ are the shape and rate parameters of the prior. Given a fixed $\boldsymbol\theta$ we will then have
$$
\pi(s | \boldsymbol\theta, {\bf g}) 
\ \propto \ \pi_{\text{exp}}({\bf g} | \boldsymbol\theta, s) \, \pi_{\text{prior}}(s) 
\ \propto \ s^N \, e^{- s \, d_W}
\, s^{a-1} \, e^{-b \, s}= s^{a+N-1} \, e^{-s (b + d_W)},
$$
where $d_W=d_W({\bf f}(\boldsymbol\theta),{\bf g})$ is the Wasserstein
metric given by \eqref{W2_distance}. Hence we obtain a Gamma posterior
on $s$:
$$
s \sim \text{Gamma}(a^*,b^*), \qquad a^*=a + N, \quad b^*=b+d_W.
$$
This will be used to generate a new sample of $s$, given a fixed
$\boldsymbol\theta$.

{\color{black}The Gamma prior used here is a conjugate prior, giving a
  Gamma posterior for $s$. We notice
  that other priors may be used for the
rate parameter $s$. Here, we use a Gamma prior, since the rate
parameter is related to the precision parameter in a Gaussian
likelihood, for which a Gamma prior is commonly used in practice. Of
course, the change of prior of $s$ will change the posterior of $s$.}

\medskip
\noindent
{\bf Metropolis-Hastings sampler.} 
The Metropolis-Hastings sampler \cite{Metropolis:53,Hastings:70} employs the proposed exponential
likelihood \eqref{W2_likelihood} with a fixed rate parameter $s$
obtained by the Gibbs sampler. The posterior from which a sequence of
samples are to be generated reads
$$
\pi(\boldsymbol\theta | s, {\bf g}) 
\ \propto \ \pi_{\text{exp}}({\bf g} | \boldsymbol\theta, s) \, \pi_{\text{prior}}(\boldsymbol\theta).
$$
Given a sample value $\boldsymbol\theta^{(i)}$, a new sample
$\boldsymbol\theta^{(i+1)}$ is generated as follows. A candidate
sample, say $\tilde{\boldsymbol\theta}$, is first generated by a proposal distribution
$q(\boldsymbol\theta^{(i)},\tilde{\boldsymbol\theta})$ from the current
sample $\boldsymbol\theta^{(i)}$. This candidate sample is then
accepted or rejected with probability (see e.g. \cite{Chib_Greenberg:95})
$$
\alpha =
\frac{\pi(\tilde{\boldsymbol\theta} | s, {\bf
    g}) \
  q(\boldsymbol\theta^{(i)},\tilde{\boldsymbol\theta})}{\pi(\boldsymbol\theta^{(i)}
  | s, {\bf g}) \ q(\tilde{\boldsymbol\theta}, \boldsymbol\theta^{(i)})} 
= \frac{
\pi_{\text{exp}}({\bf g} | \tilde{\boldsymbol\theta},s) \
\pi_{\text{prior}}(\tilde{\boldsymbol\theta}) \ q(\boldsymbol\theta^{(i)},\tilde{\boldsymbol\theta})}{\pi_{\text{exp}}({\bf g} |
\boldsymbol\theta^{(i)},s) \ \pi_{\text{prior}}(\boldsymbol\theta^{(i)}) \ q(\tilde{\boldsymbol\theta}, \boldsymbol\theta^{(i)})}.
$$
It is to be noted that the choice of prior would depend on the
problem at hand and expert's opinion. For instance, the prior may be
informative with uniform or other types of distributions, or it may be
non-informative, i.e. $\pi_{\text{prior}}(\boldsymbol\theta)=1$. We may even consider
hierarchical priors that include other hyper-priors. For the sake of
simplicity and without imposing any restriction, in the numerical
examples in Section \ref{sec:num1} we consider uniform priors on $\boldsymbol\theta$,
\begin{equation}\label{uniform_prior}
\pi_{\text{prior}}(\boldsymbol\theta) = \prod_{i=1}^m \pi_{\text{prior}}(\theta_i), \qquad \theta_i
\sim \text{Unif}(\alpha_i,\beta_i),
\end{equation}
where the bounds $\{ (\alpha_i,\beta_i) \}_{i=1}^m$ are to be selected depending on the
problem at hand. 
We further use a Gaussian random walk proposal to generate
a new sample $\tilde{\boldsymbol\theta}$ from a current sample $\boldsymbol\theta^{(i)}$:
\begin{equation}\label{proposal_random_walk}
\tilde{\boldsymbol\theta} \sim \text{Normal} (\boldsymbol\theta^{(i)}, \,
\Sigma),
\end{equation}
where $\Sigma \in {\mathbb R}^{m \times m}$ is a covariance matrix
to be selected so that the proposal distribution is neither too wide
nor too narrow. The former would result in an acceptance rate close to
zero, and the chain would rarely move to a different
sample. The latter would result in an acceptance rate close to one,
however, the generated samples would not cover the whole support of
posterior. 
We note that since such a proposal is symmetric,
i.e. $q(\boldsymbol\theta^{(i)},\tilde{\boldsymbol\theta}) =
q(\tilde{\boldsymbol\theta},\boldsymbol\theta^{(i)})$, 
the ratio $\alpha$ simplifies to $\alpha = \pi(\tilde{\boldsymbol\theta} | s,
{\bf g}) / \pi(\boldsymbol\theta^{(i)} | s, {\bf g})$. 
This ratio was used in the original version of the Metropolis
algorithm \cite{Metropolis:53}. 

\begin{algorithm}[tbh]
\caption{{\fontsize{11}{12}\selectfont
    Metropolis-Hastings-within-Gibbs Sampling in the Wasserstein-Bayesian framework}}
\label{ALG_MG}
\begin{algorithmic} 
\medskip
\STATE {\bf 1.} {\it Initialization}: Select an arbitrary point
$(\boldsymbol\theta^{(0)}, \, s^{(0)})$, and set $i=0$.

\medskip
\STATE {\bf 2.} {\it Gibbs sampler}: Generate $s^{(i+1)}$ from the
posterior $\pi (s |
\boldsymbol\theta^{(i)}, {\bf g})$ with Gamma distribution, 
$$
s^{(i+1)}  \sim \text{Gamma}(a^*,b^*), \qquad a^*=a + N, \quad
b^*=b +d_W({\bf f}(\boldsymbol\theta^{(i)}),{\bf g}).
$$

\medskip
\STATE {\bf 3.} {\it Metropolis-Hastings sampler}: Generate $\boldsymbol\theta^{(i+1)}$
as follows:

\begin{itemize}

\item Sample a candidate $\tilde{\boldsymbol\theta}$ from
a proposal distribution
$q(\boldsymbol\theta^{(i)},\tilde{\boldsymbol\theta})$.

\item Compute the ratio 
$$
\alpha(\boldsymbol\theta^{(i)},\tilde{\boldsymbol\theta}) =
\frac{
\pi_{\text{exp}}({\bf g} | \tilde{\boldsymbol\theta},s^{(i+1)}) \
\pi_{\text{prior}}(\tilde{\boldsymbol\theta}) \ q(\boldsymbol\theta^{(i)},\tilde{\boldsymbol\theta})}{\pi_{\text{exp}}({\bf g} |
\boldsymbol\theta^{(i)},s^{(i+1)}) \ \pi_{\text{prior}}(\boldsymbol\theta^{(i)}) \ q(\tilde{\boldsymbol\theta}, \boldsymbol\theta^{(i)})},
$$
with the likelihood $\pi_{\text{exp}}$ and the prior
$\pi_{\text{prior}}$ given in \eqref{W2_likelihood} and
\eqref{uniform_prior}, respectively.

\item Set 
$$
\boldsymbol\theta^{(i+1)} = \left\{ \begin{array}{l l}
\tilde{\boldsymbol\theta} & \qquad \text{if} \ \  \text{Unif}(0,1) \le \alpha(\boldsymbol\theta^{(i)},\tilde{\boldsymbol\theta}), \\
\boldsymbol\theta^{(i)} &  \qquad \text{otherwise}.
\end{array} \right.
$$
\end{itemize}
 
\medskip
\STATE {\bf 4.} {\it Iteration}: Increment $i$ by 1 and go to step {\bf 2}. 

\end{algorithmic}
\end{algorithm}

\section{Application to seismic inversion} 
\label{sec:num1} 
 
In this section, we apply the proposed Bayesian framework to the
inverse uncertainty quantification of wave propagation problems. More
specifically, we consider two problems: source inversion and
material inversion. We present several numerical examples and 
%
%
demonstrate the applicability and performance of the proposed
framework in comparison with the conventional Bayesian framework
based on Gaussian likelihoods.

\subsection{A one-dimensional source inversion problem}

{\color{black}This example concerns a simple one-dimensional source inversion
problem where an initial wave pulse propagates with a constant
speed. 
The purpose of this example is to demonstrate the two major advantages of the
Wasserstein distance over the standard least-squares norm separately, i.e. 1)
handling complex noise structures, and 2) handling the problem with
phase shifts and the cycle skipping issue. Correspondingly, we split the example into two parts: A) with a complicated noise
structure and unknown amplitude, and B) with additive Gaussian noise
and unknown phase and amplitude. Since the only parameter in part A is the
amplitude, there will be no phase shift and no cycle skipping
issue. Therefore, part A is designed to demonstrate the advantage of the Wasserstein
distance only in handling complicated noise structures. In part B, on the
other hand, noise is not an issue. With phase being an unknown, part B
focuses on the performance of the two distances in handling the cycle skipping
issue and the presence of multiple local minima.

\medskip
\noindent
{\bf Problem formulation.} 
Consider the Cauchy problem for the 1D wave equation 
\begin{align}
&u_{tt}(t,x) - u_{xx}(t,x) = 0, \qquad t \in [0,T], \quad x \in {\mathbb
R},\\
&u(0,x)=h(x; x_0, a), \qquad u_t(0,x)=0.
\end{align}
The initial data
$$
h(x; x_0, a) = a \, \left( e^{-100 (x-x_0-0.5)^2} +
  e^{-100 (x-x_0)^2} + e^{-100 (x-x_0+0.5)^2} \right), 
$$
acts as a source term generating an initial wave pulse, involving two
parameters: the initial location $x_0$ and the amplitude $a$ of the
wave pulse. 
%
The solution to this simple problem is given by d'Alembert's formula
\begin{equation}\label{dalembert}
u(t,x; x_0, a)= \frac{1}{2} \, h(x-t; x_0, a) + \frac{1}{2} \, h(x+t; x_0, a). 
\end{equation}

\subsubsection{A complicated noise structure and unknown amplitude}

\noindent
{\bf Synthetic data.} Let $x_0 = 0$ and $\theta = a$ be the only unknown parameter. Consider a uniform array of $N_r = 7$ receivers located at
$$
x_1=-3, \ \ x_2=-2, \ \ x_3=-1, \ \ x_4=0, \ \ x_5=1, \ \ x_6=2,
\ \ x_7=3.
$$ 
Each receiver located at $x_r$, with $r=1, \dotsc, N_r$, records a
noisy discrete-time signal $g(t_k,x_r)$ over the time interval $[0,T]$
at $N$ discrete time levels $t_k = (k-1) \, \Delta t$, with $\Delta t = T/(N - 1)$ and
$k=1, \dotsc, N$. Let $f(t_k,x_r; \theta)$ denote the
corresponding simulated signal for a given $\theta$, computed by
\eqref{dalembert} with $x_0=0$. We select a fixed parameter 
$$
\theta^*=5,
$$
and generate synthetic data as follows. We first compute $f(t_k,x_r; \theta^*)$ by \eqref{dalembert}
and then generate $g(t_k, x_r)$ by polluting $f(t_k,x_r; \theta^*)$ with a noise of the form \eqref{noise_structure}: we multiply $f$ by
i.i.d. Gamma noise and then add an i.i.d uniform noise to the product,
that is,
$$
g(t_k, x_r) = \varepsilon_{kr}^{(1)}  \, f(t_k, x_r; \theta^*) +
\varepsilon_{kr}^{(2)}, \qquad \varepsilon_{kr}^{(1)} \sim
\text{Gamma}(60, 60), \qquad \varepsilon_{kr}^{(2)} \sim
\text{Unif}(-0.25, 0.25).
$$
Let ${\bf g}_r \in {\mathbb R}^{N}$ denote the discrete-time signal recorded at $x_r$, $r=1, \dotsc, N_r$. We have a set
of $N_r=7$ such signals, shown in Figure \ref{observations} with $T=5$
and $N = 101$. 
Similarly, we let ${\bf f}_r(\theta) \in {\mathbb R}^{N}$
be the $r$-th simulated signal, consisting of $f(t_k,x_r; \theta)$
at all $N$ discrete time levels $\{ t_k \}_{k=1}^{N}$. 
\begin{figure}[!h]
\vspace{-.2cm}
\begin{center}
\includegraphics[width=0.6\linewidth]{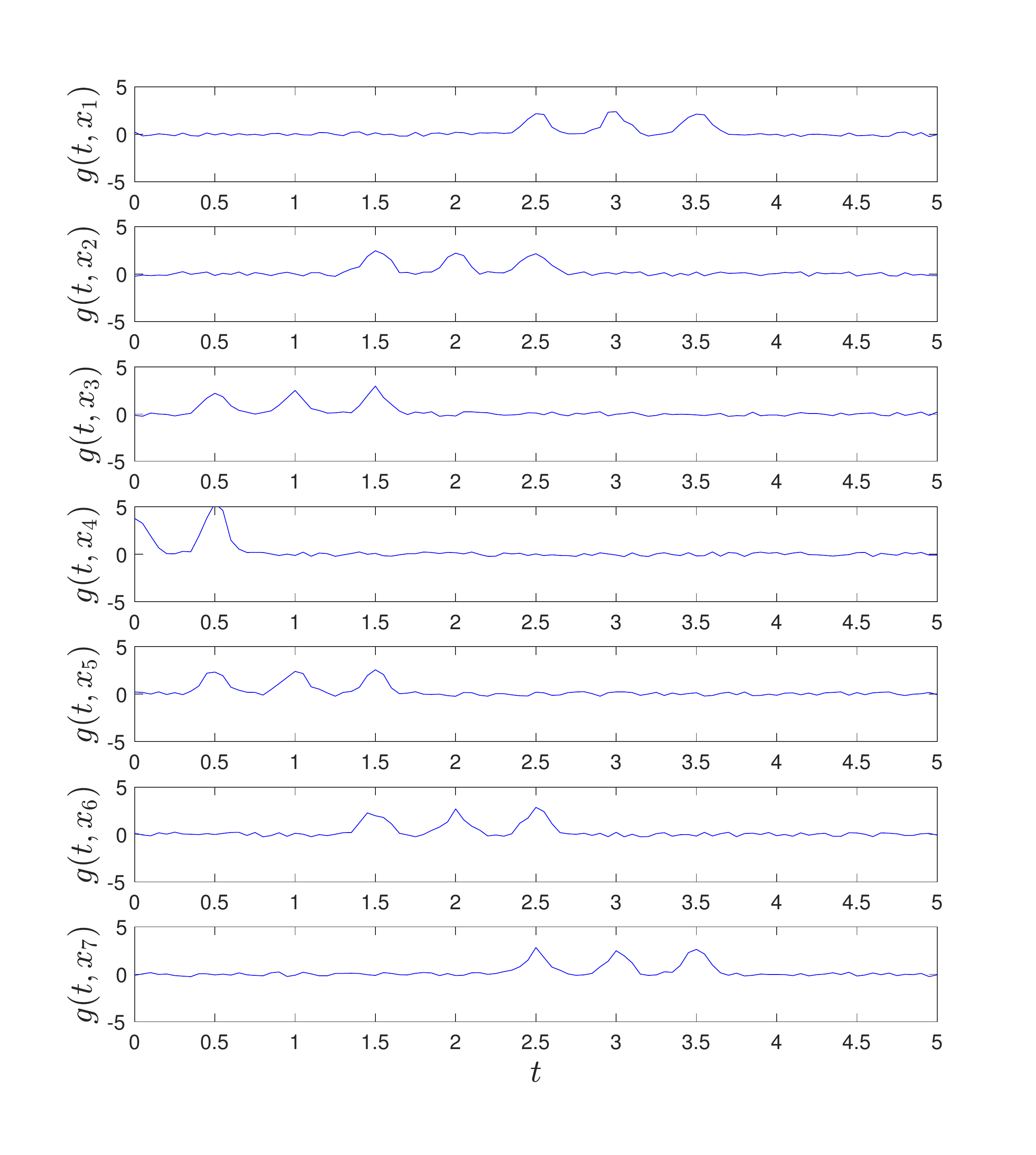}     
\vspace{-1cm}
\caption{Observed signals with complex noise structures, recorded at seven receivers.}
\label{observations}
\end{center}
\end{figure}
%
%



\medskip
\noindent
{\bf Computations.} The goal is to employ Bayesian inversion and
compute the conditional posterior of the amplitude parameter
$\theta$. 
We follow Algorithm \ref{ALG_MG} with the
following choices:
\begin{itemize}
\item {\it Likelihood}: 
$
\pi_{\text{exp}}( {\bf g} | \theta) =s^N \, \exp ( -s
\, \sum_{r=1}^{N_r} d_W({\bf f}_r(\theta), {\bf g}_r)).
$

\item {\it Priors}: 
$
\theta \sim \text{Unif}(2 ,8),\ \ \ \
s \sim \text{Gamma}(1 ,0.1).
$

\item {\it Proposal}: A Gaussian random walk
  \eqref{proposal_random_walk} with variance $\Sigma = 0.005$.

\item {\it Initial point}: $\theta^{(0)} = 3$, \ \ \
$s^{(0)}=70$.
\end{itemize}

We also consider the standard normal likelihood based on the Euclidean distance:
$$
\pi_{\text{norm}}( {\bf g} | \theta) =\frac{s^{N/2}}{(2
  \pi)^{N/2}} \, \exp ( -\frac{s}{2} 
\, \sum_{r=1}^{N_r} d_2({\bf f}_r(\theta), {\bf g}_r)),
\qquad s = \sigma^{-2},
$$
where $d_2({\bf f}_r(\theta), {\bf g}_r) = \sum_{k=1}^{N}
| g(t_k,x_r) - f(t_k,x_r;\theta)|^2$. 
We note that in this case $s = \sigma^{-2}$ will be the precision parameter.

We run the algorithm with $M_0=30000$ iterations and remove the first
$M_b=10000$ samples, known as the burn-in period. We also use a
thinning period of $M_t=4$, that is, we keep every 4th samples and
discard the rest. This would give a total of $M=5000$ Markov chain samples. 
Figure \ref{fig_ex1_results} shows the posterior histograms and the
trace plots of the Markov chain samples generated by the algorithm
with both (a) the proposed likelihood and (b) the Gaussian likelihood. We observe
that while both posteriors are centered around the true noise-free
parameter, the one obtained by the proposed exponential
likelihood function is much more concentrated compared to the one
obtained by the Gaussian likelihood function. This clearly demonstrate
the power of the proposed likelihood in capturing more information in
the available data compared to standard Gaussian likelihood. 
\begin{figure}[!h]
\center
\subfigure[Wasserstein-driven exponential likelihood.]{\includegraphics[width=0.45\linewidth]{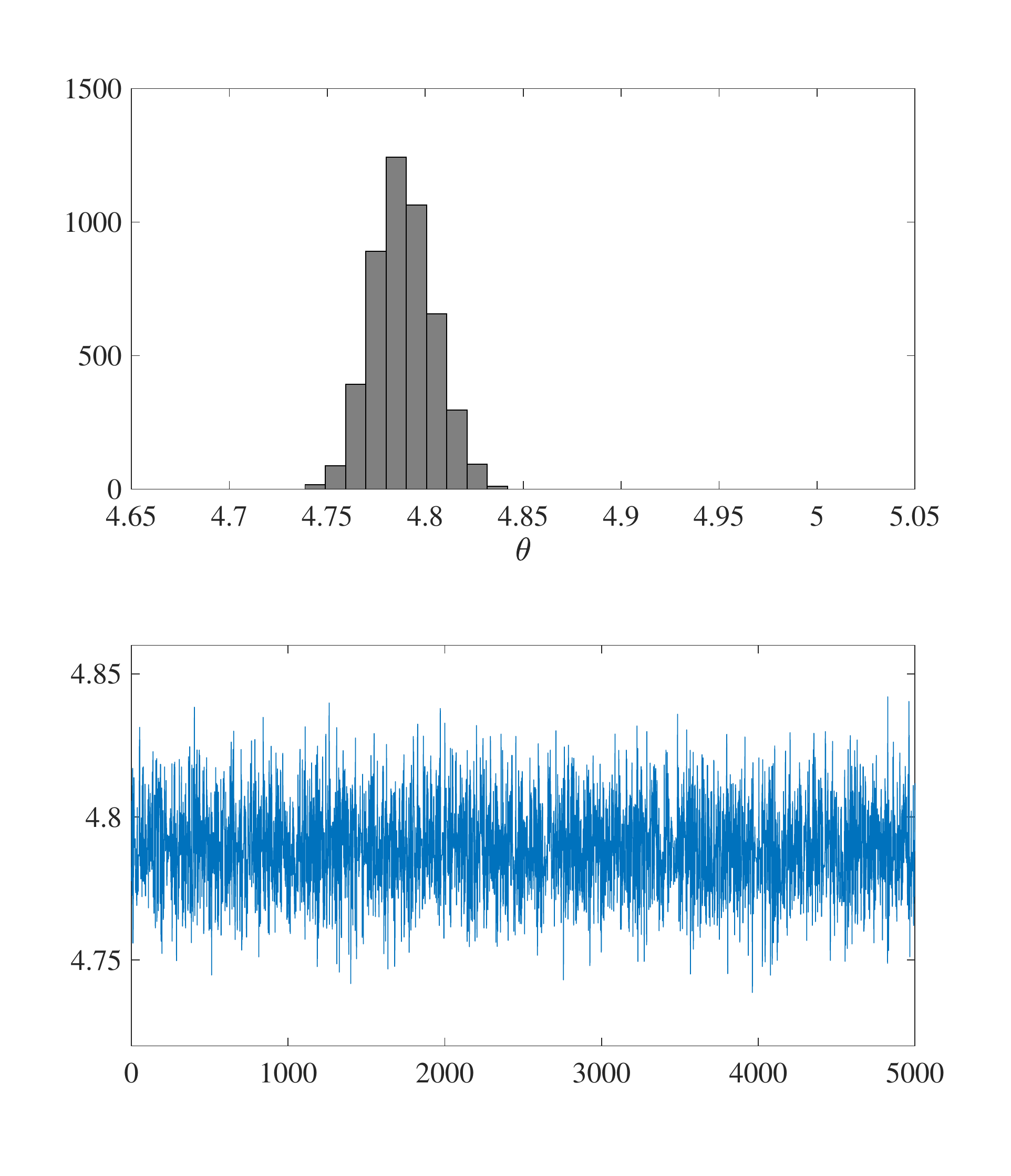}}
\subfigure[$L_2$-driven Gaussian likelihood.]{\includegraphics[width=0.45\linewidth]{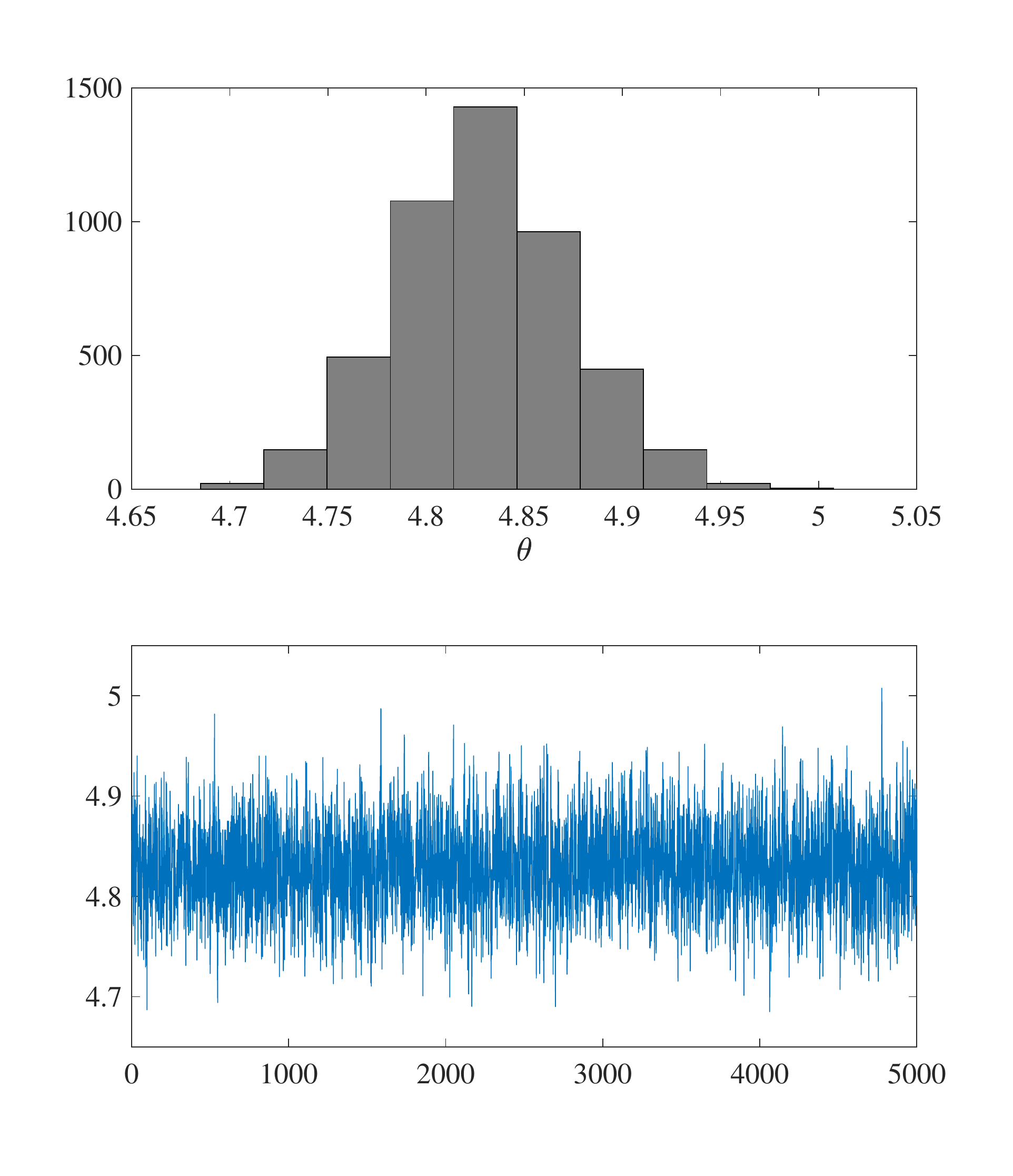}}
\caption{Posterior histograms and the trace plots of the Markov chain
  samples in two cases: (a) with the proposed
  Wasserstein-driven exponential likelihood, and (b) with the standard
  Gaussian likelihood. Clearly, the posterior obtained by the Gaussian
  likelihood function is much more concentrated compared to the one
obtained by the Gaussian likelihood function.} 
\label{fig_ex1_results}
\end{figure}

\subsubsection{An additive Gaussian noise and unknown phase and amplitude}

\noindent
{\bf Synthetic data.} Let $\boldsymbol\theta = (x_0, a)$ be the unknown
parameter vector. Consider a uniform array of $N_r = 7$ receivers located at
$$
x_1=-3, \ \ x_2=-2, \ \ x_3=-1, \ \ x_4=0, \ \ x_5=1, \ \ x_6=2,
\ \ x_7=3.
$$ 
Each receiver located at $x_r$, with $r=1, \dotsc, N_r$, records a
noisy discrete-time signal $g(t_k,x_r)$ over the time interval $[0,T]$
at $N$ discrete time levels $t_k = (k-1) \, \Delta t$, with $\Delta t = T/(N - 1)$ and
$k=1, \dotsc, N$. Let $f(t_k,x_r; \boldsymbol\theta)$ denote the
corresponding simulated signal for a given $\boldsymbol\theta$, computed by \eqref{dalembert}. We select a fixed parameter vector 
$$
\boldsymbol\theta^*=(0,5),
$$
and generate synthetic data as follows. We first compute $f(t_k,x_r; \boldsymbol\theta^*)$ by \eqref{dalembert}
and then generate $g(t_k, x_r)$ by polluting $f(t_k,x_r;
\boldsymbol\theta^*)$ with an additive Gaussian noise,
$$
g(t_k, x_r) = f(t_k, x_r; \boldsymbol\theta^*) +
\varepsilon_{kr}, \qquad \varepsilon_{kr} \sim
\text{Normal}(0, 0.1).
$$
Let ${\bf g}_r \in {\mathbb R}^{N}$ denote the discrete-time signal recorded at $x_r$, $r=1, \dotsc, N_r$. We have a set
of $N_r=7$ such signals, shown in Figure \ref{observations} with $T=5$
and $N = 101$. 
Similarly, we let ${\bf f}_r(\boldsymbol\theta) \in {\mathbb R}^{N}$
be the $r$-th simulated signal, consisting of $f(t_k,x_r; \boldsymbol\theta)$
at all $N$ discrete time levels $\{ t_k \}_{k=1}^{N}$. 
\begin{figure}[!h]
\vspace{-.2cm}
\begin{center}
\includegraphics[width=0.6\linewidth]{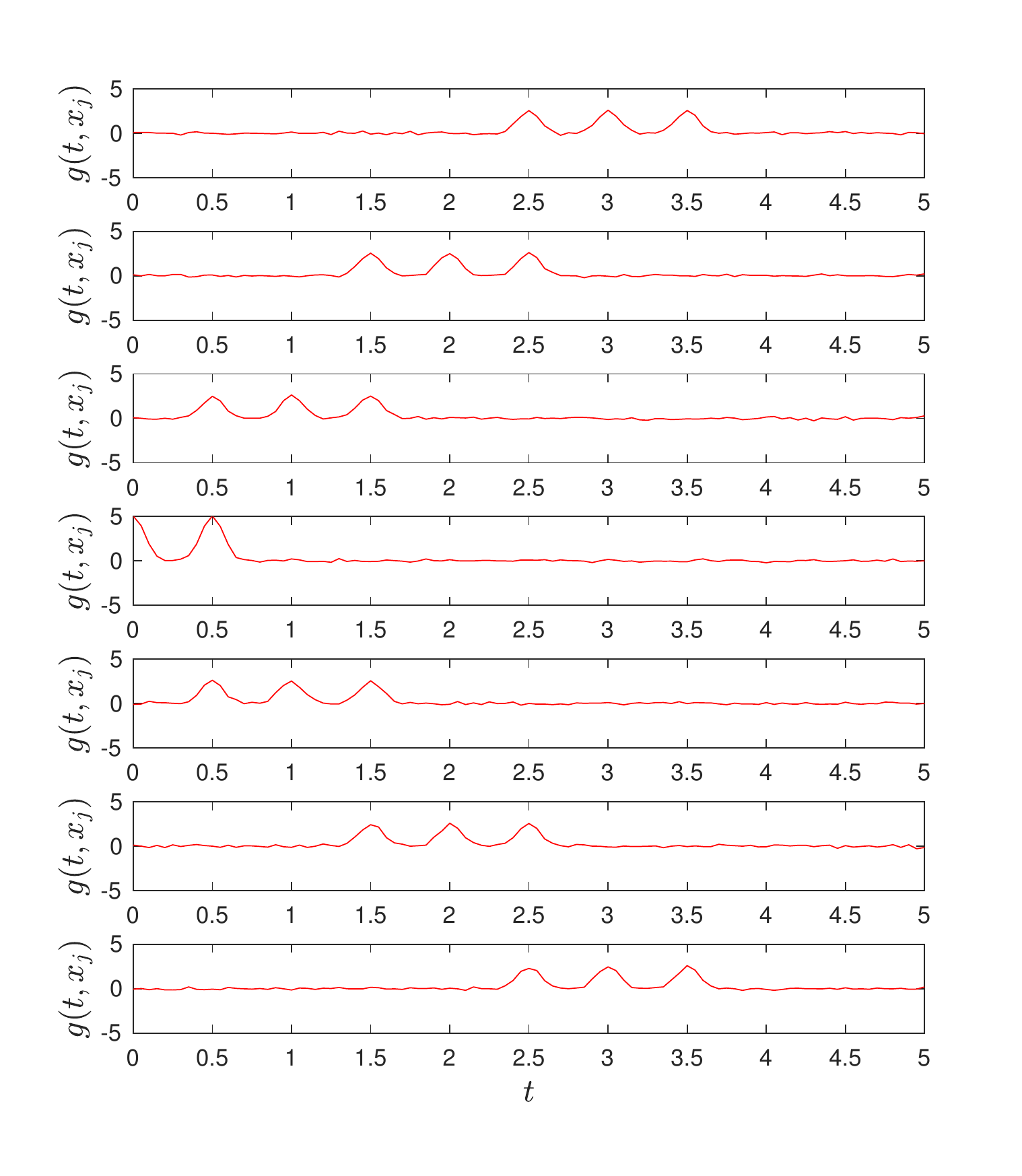}     
\vspace{-1cm}
\caption{Observed signals with additive Gaussian noise, recorded at seven receivers.}
\label{observations}
\end{center}
\end{figure}
%

\noindent
{\bf Computations.} The goal is to employ Bayesian inversion and
compute the conditional posterior of the model parameters
$\boldsymbol\theta$. 
We follow Algorithm \ref{ALG_MG} with the
following choices:
\begin{itemize}
\item {\it Likelihood}: 
$
\pi_{\text{exp}}( {\bf g} | \boldsymbol\theta) =s^N \, \exp ( -s
\, \sum_{r=1}^{N_r} d_W({\bf f}_r(\boldsymbol\theta), {\bf g}_r)).
$

\item {\it Priors}: 
$
\theta_1 \sim \text{Unif}(-3 ,3), \ \ \ \
\theta_2 \sim \text{Unif}(2 ,8),\ \ \ \
s \sim \text{Gamma}(1 ,0.1).
$

\item {\it Proposal}: A Gaussian random walk
  \eqref{proposal_random_walk} with covariance $\Sigma = \text{diag}(0.005, 0.005)$.

\item {\it Initial point}: $\boldsymbol\theta^{(0)} = (0.6, 3)$, \ \ \
$s^{(0)}=70$.
\end{itemize}

We also consider the standard normal likelihood based on the Euclidean distance:
$$
\pi_{\text{norm}}( {\bf g} | \boldsymbol\theta) =\frac{s^{N/2}}{(2
  \pi)^{N/2}} \, \exp ( -\frac{s}{2} 
\, \sum_{r=1}^{N_r} d_2({\bf f}_r(\boldsymbol\theta), {\bf g}_r)),
\qquad s = \sigma^{-2},
$$
where $d_2({\bf f}_r(\boldsymbol\theta), {\bf g}_r) = \sum_{k=1}^{N}
| g(t_k,x_r) - f(t_k,x_r;\boldsymbol\theta)|^2$. 
We note that in this case $s = \sigma^{-2}$ will be the precision parameter.

We run the algorithm with $M_0=25000$ iterations and remove the first
$M_b=5000$ samples, known as the burn-in period. We also use a
thinning period of $M_t=4$, that is, we keep every 4th samples and
discard the rest. This would give a total of $M=5000$ Markov chain samples. 
Figure \ref{fig_ex1_results} shows the posterior histograms and the
trace plots of the Markov chain samples generated by the algorithm
with both (a) the proposed likelihood and (b) the Gaussian likelihood. We observe
that the samples obtained by the standard Gaussian likelihood converge
to wrong distributions,
while the Markov chain samples obtained by the proposed exponential
likelihood function are centered around the true, noise-free parameters. 
\begin{figure}[!ht]
\center
\subfigure[Wasserstein-driven exponential likelihood.]{\includegraphics[width=0.65\linewidth]{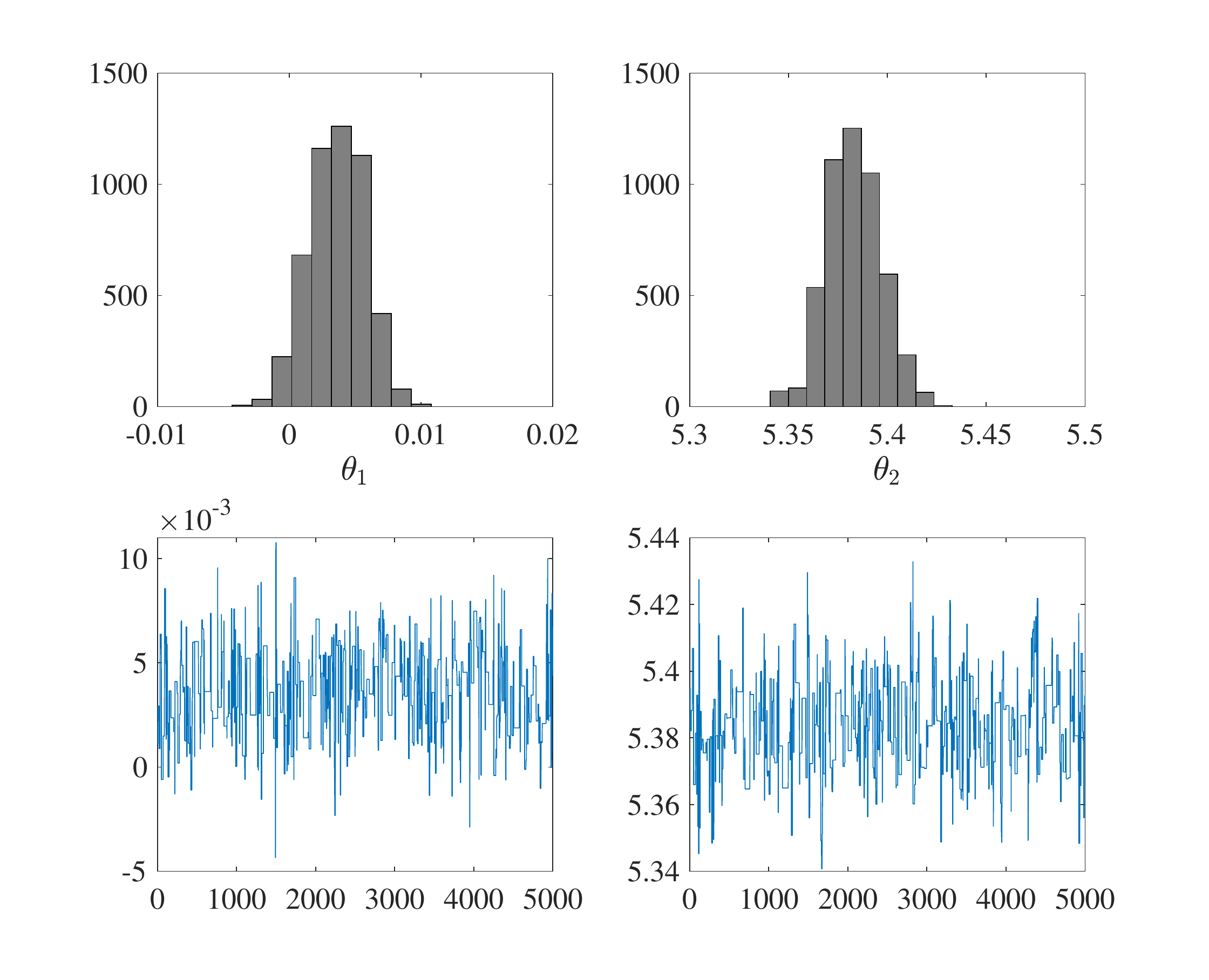}}
\subfigure[$L_2$-driven Gaussian likelihood.]{\includegraphics[width=0.65\linewidth]{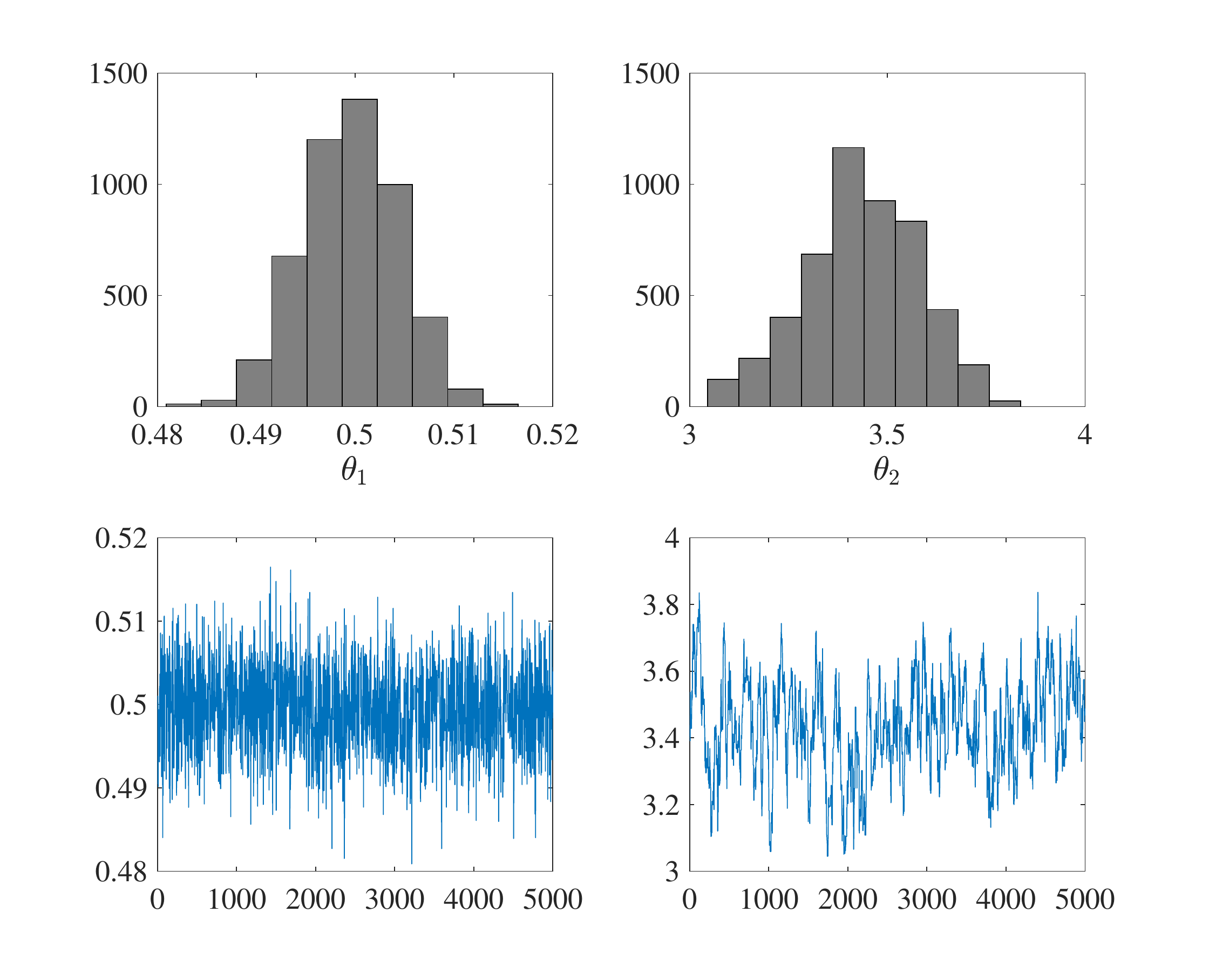}}
\caption{Posterior histograms and the trace plots of the Markov chain
  samples in two cases: (a) with the proposed
  Wasserstein-driven exponential likelihood, and (b) with the standard
  Gaussian likelihood. Clearly, the samples obtained by the Gaussian
  likelihood function do not converge to the correct distribution,
  while the samples obtained by the proposed likelihood are accumulated
  close to the true, noise-free parameters.} 
\label{fig_ex1_results}
\end{figure}

The wrong convergence of the standard Bayesian approach, as shown in Figure
\ref{fig_ex1_results}, is due to the
existence of multiple local maxima in the Gaussian likelihood
function. Consequently, the Markov chain samples may get trapped in a
local maximum and hence miss the global maximum. 
To illustrate this, in Figure \ref{fig_likelihoods} we plot the log-likelihood functions
$\log \, L_{exp}(\theta_1,5)$ and $\log \, L_{norm}(\theta_1,5)$ versus $\theta_1 \in
[-3 ,3]$. We observe that the Gaussian likelihood has multiple local
maxima, while the proposed exponential likelihood has a more-convex
shape, making it easy for MCMC to infer its global maximum.
\begin{figure}[!ht]
\center
\subfigure[Wasserstein-driven exponential
likelihood.]{\includegraphics[width=0.49\linewidth]{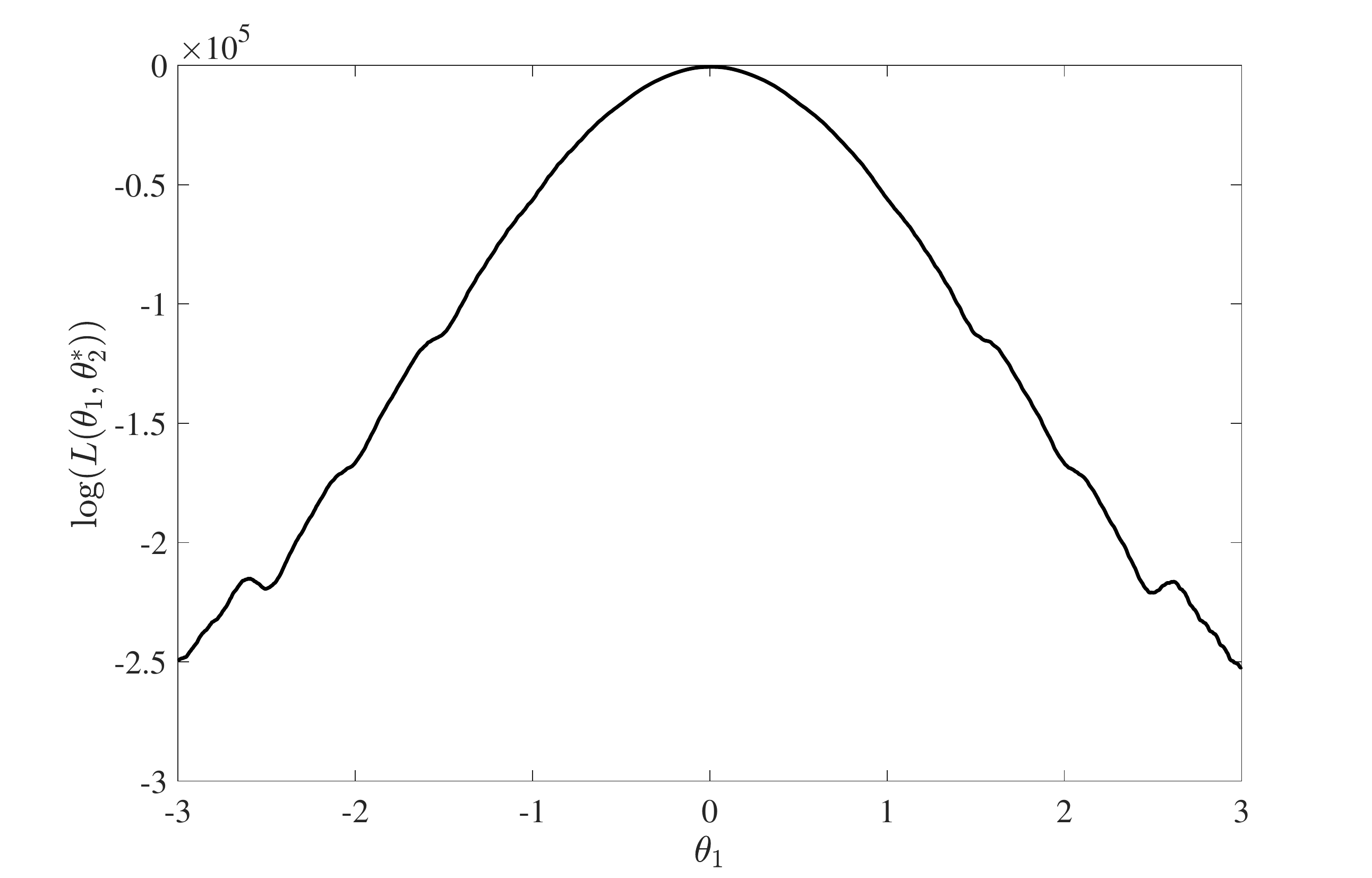}}
\hskip .01cm
\subfigure[$L_2$-driven Gaussian likelihood.]{\includegraphics[width=0.49\linewidth]{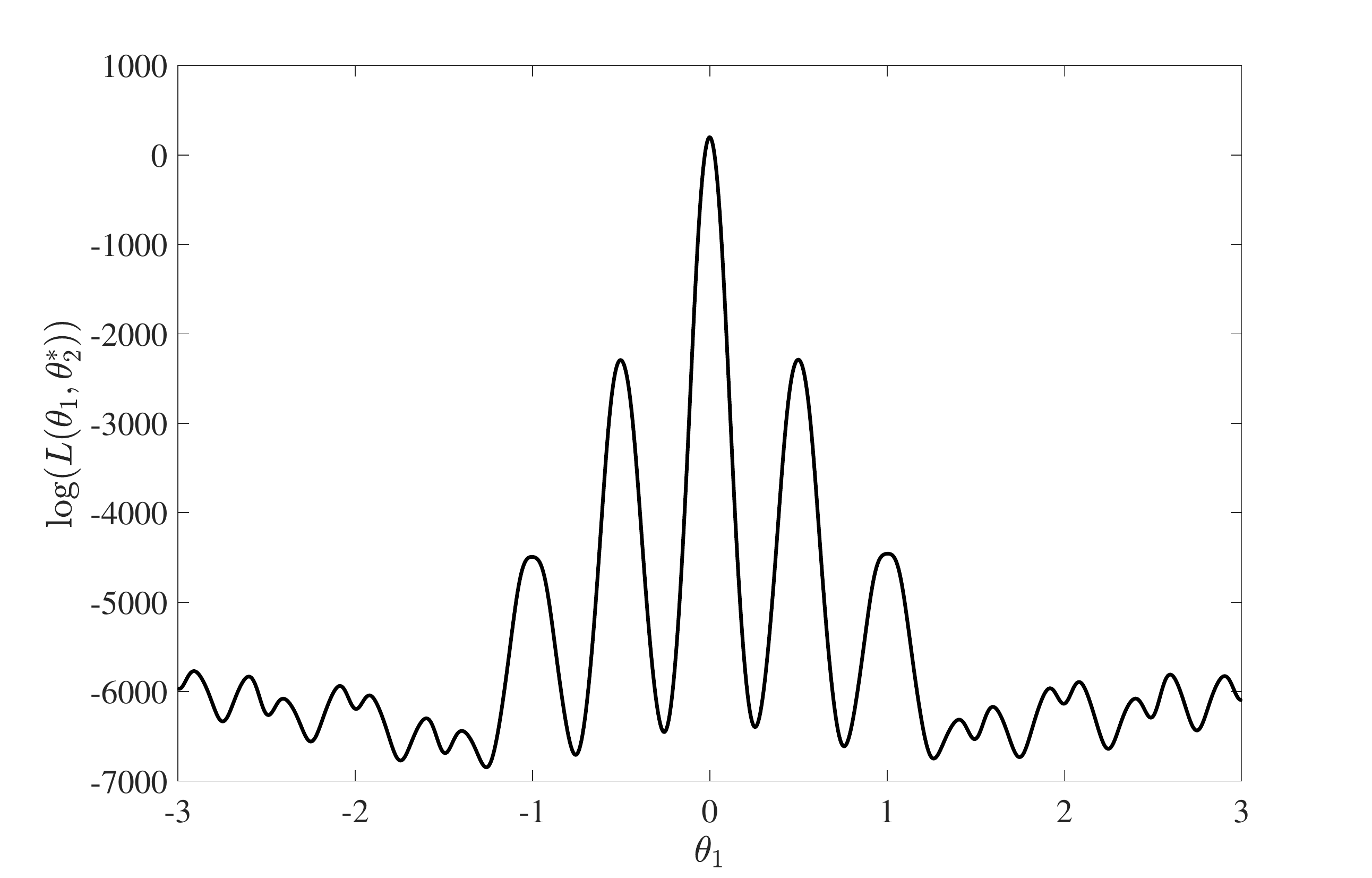}}
\caption{The log-likelihood functions $\log \, L_{exp}(\theta_1,5)$ and
  $\log \, L_{norm}(\theta_1,5)$ versus $\theta_1 \in [-3 ,3]$. 
} 
\label{fig_likelihoods}
\end{figure}


}
\subsection{A two-dimensional seismic source inversion problem}

This example concerns the inference of a seismic source in a layered
isotropic elastic material in two dimensions. The source models a
simplified earthquake, represented by four unknown parameters.

\medskip
\noindent
{\bf Problem formulation.} 
We consider the initial-boundary value problem (IBVP) for the elastodynamic
wave equation, 
\begin{subequations}\label{elastic}
\begin{align}
&\nu(\bold{x}) \, {\bf u}_{tt} (t,\bold{x})- \nabla \cdot \boldsymbol \sigma ({\bf u}(t,\bold{x})) = {\bf h}(t,\bold{x};\boldsymbol{\theta}) \hskip 1cm  \text{in  } [0,T] \times D, \label{elastic_1st}\\
&{\bf u}(0,\bold{x})={\bf 0}, \hskip 6mm {\bf u}_t(0,\bold{x})={\bf 0} \hskip 3.37cm \text{on } \{ t=0 \} \times D, \label{elastic_2nd}\\
&\boldsymbol \sigma({\bf u}(t,\bold{x})) \cdot \hat{\bf n} = {\bf 0} \hskip 5.16cm \text{on } [0,T] \times \partial D_0, \label{elastic_3rd}\\
&{\bf u}_t (t,\bold{x}) = B({\bf x}) \, \boldsymbol \sigma({\bf u}(t,\bold{x})) \cdot \hat{\bf n} \hskip 3.12cm \text{on } [0,T] \times \partial D_1. \label{elastic_4th}
\end{align}
\end{subequations}
The solution ${\bf u}= (u_1, u_2)^{\top}$ represents the
displacement field, $t \in [0,T]$ is the time, ${\bf x}= (x_1, x_2)
\in D$ is the location, and $\boldsymbol \sigma$ is the stress tensor given by
\begin{equation}\label{stress}
 \boldsymbol \sigma({\bf u}) = \lambda (\bold{x}) \, \nabla \cdot {\bf u} \, I+ \mu (\bold{x}) \, 
 (\nabla {\bf u} + (\nabla {\bf u})^{\top}),
\end{equation}
where $I$ is the identity matrix. The computational domain is a box $D = [-10000,10000] \times [-15000,0]$, consisting of two
layers: the top layer $D_{I}$ extends over $-1000 \le x_2 \le 0$, and the bottom layer $D_{II}$ is given by $x_2 \le -1000$. 
The material properties are characterized by the density $\nu$ and the
Lam{\'e} parameters, $\lambda$ and $\mu$. 
The system 
\eqref{elastic} admits longitudinal (or pressure) and transverse (or shear) 
waves, which, in the case of constant density, propagate at the wave speeds
$$
c_p = \sqrt{(2 \, \mu+\lambda) / \nu}, \qquad c_s = \sqrt{\mu / \nu},
$$
respectively. 
The material density and velocities are assumed to be known, given by  
\begin{eqnarray*}
\nu({\bf x})=  \left\{ \begin{array}{l l}
2600 &  {\bf x} \in D_I\\
2700 & {\bf x} \in D_{II}
\end{array} \right.,
\quad
c_p({\bf x})=  \left\{ \begin{array}{l l}
4000 &  {\bf x} \in D_I\\
6000 & {\bf x} \in D_{II}
\end{array} \right.,
\quad
c_s({\bf x})=  \left\{ \begin{array}{l l}
2000 &  {\bf x} \in D_I\\
3464 & {\bf x} \in D_{II}
\end{array} \right..
\end{eqnarray*}
We impose a homogeneous Neumann (stress-free) boundary condition
\eqref{elastic_3rd} on the free surface $\partial D_0 =\{ {\bf x} |x_2
=0 \}$, and the first-order Clayton-Engquist non-reflecting boundary conditions
\cite{Clayton_Engquist:1977} on the other three boundaries $\partial
D_1 = \partial D\setminus \partial D_0$. Here, ${\hat {\bf n}}$ is the
outward unit normal to the boundary, and $B$ is a given matrix related
to the non-reflecting boundary condition.

The function ${\bf h}$ represents the seismic source. We consider the
case of a point moment tensor source that models a simplified earthquake,
\begin{equation}\label{source}
{\bf h}(t,{\bf x};\boldsymbol{\theta}) = S(t) \, M \, \nabla \delta({\bf x} -{\bf x}_s),
\end{equation}
located at ${\bf x}_s = (0, z_{s})$, where $\nabla \delta$ is the
gradient of the Dirac distribution. The source time function $S(t)$
and the moment tensor are assumed to be
$$
S(t) = \frac{\omega_s}{\sqrt{2 \, \pi}} \, e^{- \omega_s^2 \, (t -
  t_s)^2 /2}, 
\qquad
M=  m_s \, \left( \begin{array}{c c}
10^{10} &  0.8 \times 10^{10}\\
0.8 \times 10^{10} &  10^{10}
\end{array} \right).
$$
parametrized by the frequency $\omega_s$, the center time $t_s$, and a
constant $m_s$. 
Hence, the source parameter vector $\boldsymbol{\theta}$ consists of four parameters,
$$
\boldsymbol{\theta} = (z_s, t_s, \omega_s, m_s).
$$

\begin{figure}[!h]
\begin{center}
\includegraphics[width=0.4\linewidth]{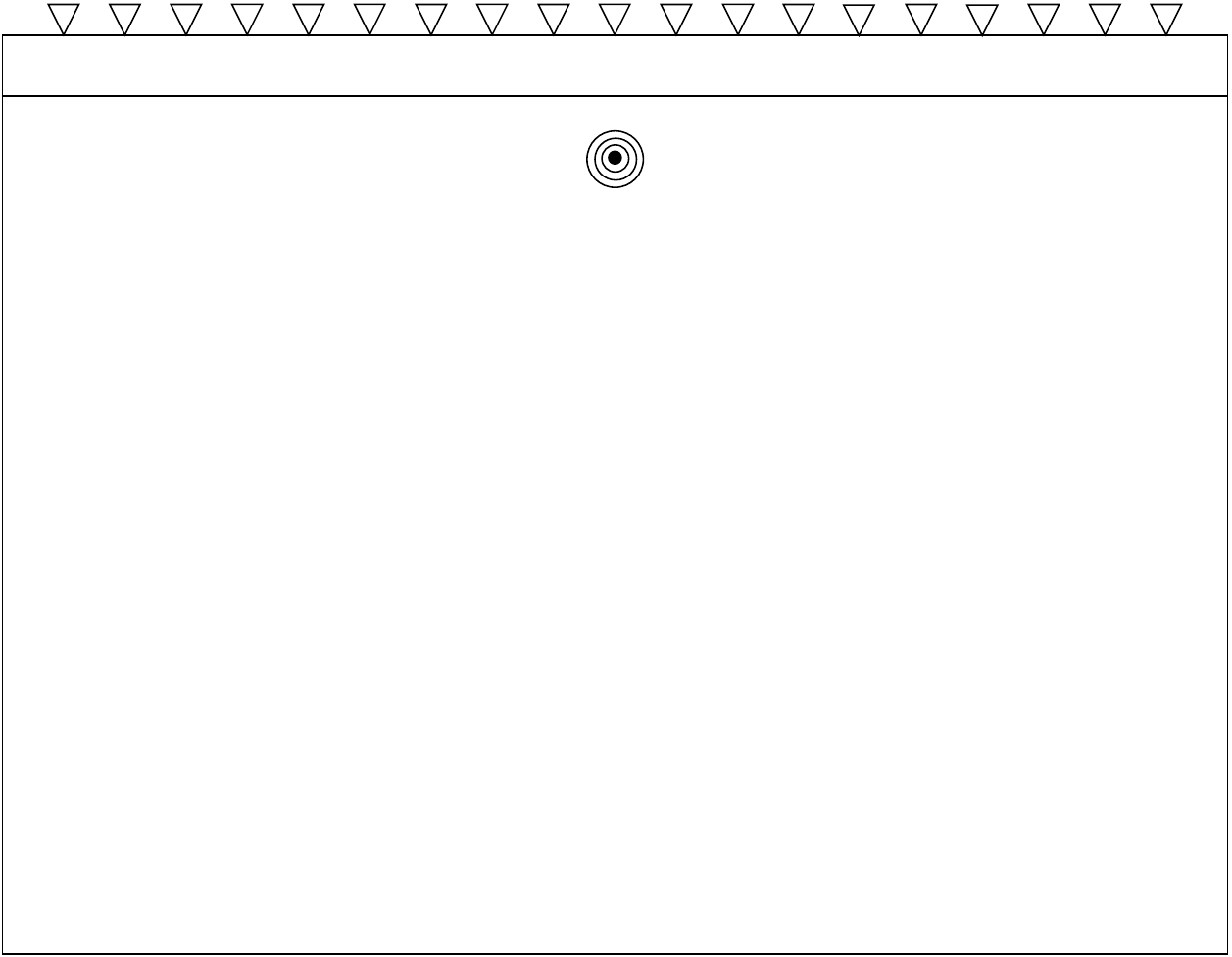}     
\vspace{-.01cm}
\caption{A schematic representation of the computational domain,
  consisting of two layers, the unknown seismic source (concentric circles), and the uniform array of receivers (black
  $\triangledown$) placed at the ground surface.}
\label{EX2_domain_setting}
\end{center}
\end{figure}

\medskip
\noindent
{\bf Synthetic data.} We place a uniform array of $N_r = 90$ receivers located at
$$
{\bf x}_r = (x_{1,r}, x_{2,r}) = (-8900 + 200 (r-1), \ 0), \qquad r=1, 2, \dotsc, N_r.
$$ 
Figure \ref{EX2_domain_setting} shows the computational domain and the
configuration of the problem. 
Each receiver located at ${\bf x}_r$, with $r=1, \dotsc, N_r$, records a
noisy vector-valued discrete-time signal $\bar{g}(t_k,{\bf x}_r) =
(g_1(t_k,{\bf x}_r), g_2(t_k,{\bf x}_r))$ over the time interval $[0,T]$
at $N$ discrete time levels $t_k = (k-1) \, \Delta t$, with $\Delta t = T/(N - 1)$ and
$k=1, \dotsc, N$. Moreover, let $\bar{f}(t_k,{\bf x}_r;
\boldsymbol\theta)= (f_1(t_k,{\bf x}_r; \boldsymbol\theta),
f_2(t_k,{\bf x}_r; \boldsymbol\theta))$ denote the
corresponding simulated signal for a given $\boldsymbol\theta$ and
computed by numerically solving \eqref{elastic}. We employ a
second-order accurate finite difference scheme, proposed in
\cite{Nilsson_etal:07}, with grid size $\Delta x_1 = \Delta x_2 = 100$
and time step $\Delta t = 0.0125$. We set the final time $T=8$ with $N
= 641$ time steps. We also employ the technique
proposed in \cite{Sjogreen_Petersson:14} for regularized
approximations of the Dirac distribution and its gradient to achieve
point-wise convergence of the solution away from the singular source. We refer to \cite{LMT:15} for details of the employed forward solver. 
We then choose a fixed parameter vector 
$$
\boldsymbol\theta^*=(-2000,1,4,10^4),
$$
and generate synthetic data as follows. We first numerically compute
$\bar{f}(t_k,{\bf x}_r; \boldsymbol\theta^*)$ and then generate
$\bar{g}(t_k, {\bf x}_r)$ by polluting $\bar{f}(t_k,{\bf x}_r; \boldsymbol\theta^*)$
with a noise of the form \eqref{noise_structure}: 
$$
g_{1,2}(t_k, {\bf x}_r) = \varepsilon_{kr}^{(1)}  \, f_{1,2}(t_k, {\bf x}_r; \boldsymbol\theta^*) +
\varepsilon_{kr}^{(2)}, \qquad \varepsilon_{kr}^{(1)} \sim
\text{Gamma}(500, 500), \qquad \varepsilon_{kr}^{(2)} \sim
\text{Unif}(-\frac1{80}, \frac1{80}).
$$
Let ${\bf g}_r \in {\mathbb R}^{N}$ denote the discrete-time signal recorded at the
receiver location ${\bf x}_r$, with $r=1, \dotsc, N_r$. We note that since
there are two (a horizontal and a vertical) components of the solution, we
have a set of $2 N_r=180$ such signals. Similarly, we let ${\bf f}_r(\boldsymbol\theta) \in {\mathbb R}^{N}$
be the $r$-th simulated signal, consisting of $f_1(t_k,{\bf x}_r;
\boldsymbol\theta)$ or $f_2(t_k,{\bf x}_r; \boldsymbol\theta)$ at all $N$ discrete time levels $\{ t_k \}_{k=1}^{N}$.

\medskip
\noindent
{\bf Computations.} The goal is to employ Bayesian inversion and
compute the conditional posterior of the model parameters
$\boldsymbol\theta$. 
We follow Algorithm \ref{ALG_MG} with the
following choices:
\begin{itemize}
\item {\it Likelihood}: 
$
\pi_{\text{exp}}( {\bf g} | \boldsymbol\theta) =s^N \, \exp ( -s
\, \sum_{r=1}^{2 N_r} d_W({\bf f}_r(\boldsymbol\theta), {\bf g}_r)).
$

\item {\it Priors}: 
$
\ \ \ \theta_1 \sim \text{Unif}(-2200,-1800), \ \ \ \
\theta_2 \sim \text{Unif}(0.5,1.5),\ \ \ \
\theta_3 \sim \text{Unif}(2,6),$

$\hskip 1.87cm \theta_4 \sim \text{Unif}(9600,10400),\ \ \ \
s \sim \text{Gamma}(1 ,0.1).
$

\item {\it Proposal}: A Gaussian random walk
  \eqref{proposal_random_walk} with covariance
  $\Sigma$, which is initially chosen to be diagonal and
  then updated based on the first 1000 MCMC samples.

\item {\it Initial point}: $\boldsymbol\theta^{(0)} =
  (-2100,0.7,3,9800)$, \ \ \
$s^{(0)}=2000$.
\end{itemize}

\begin{figure}[!h]
\center
\includegraphics[width=0.99\linewidth]{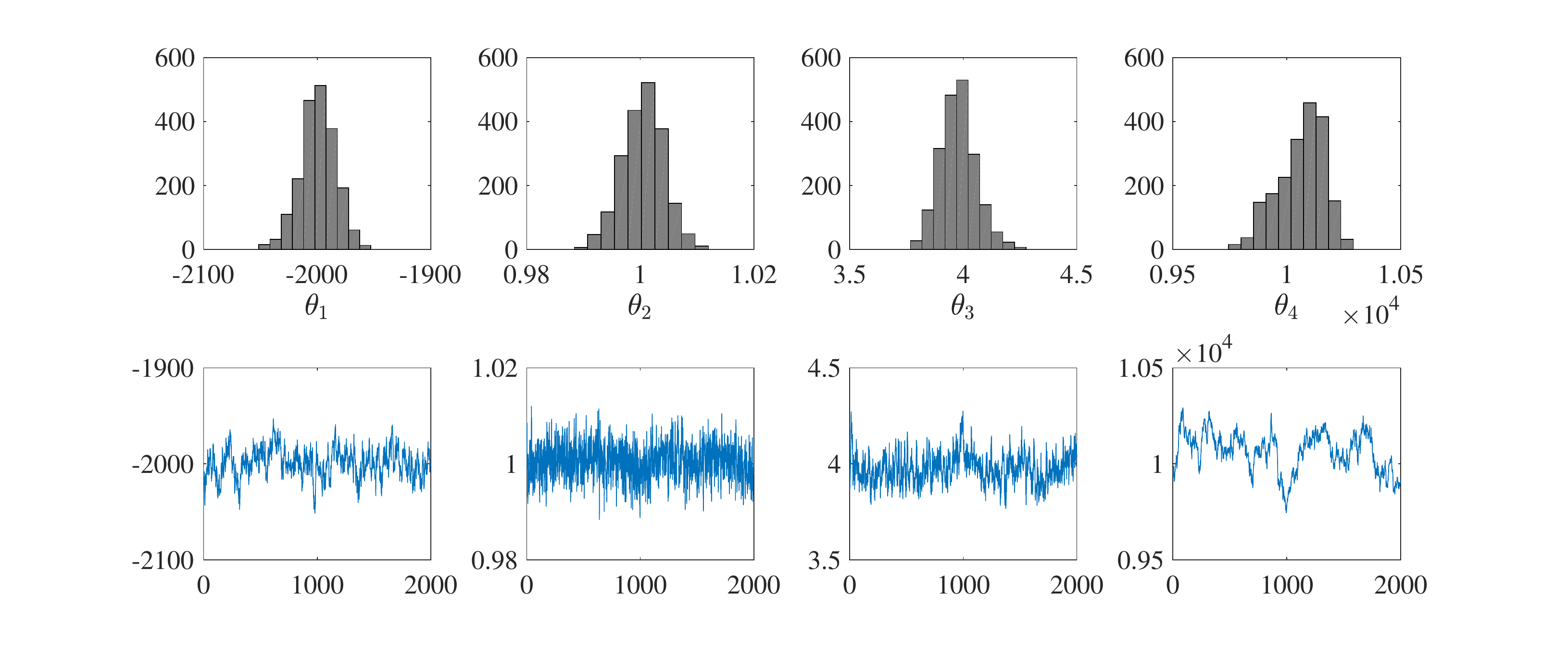}
\caption{Posterior histograms and the trace plots of the Markov chain
  samples computed by the proposed Wasserstein-driven Bayesian technique.} 
\label{fig_ex2_results}
\end{figure}

We run the algorithm with $M_0=10000$ iterations and remove the first
$M_b=2000$ samples. We also use a
thinning period of $M_t=4$. This would give a total of $M=2000$ Markov chain samples. Figure \ref{fig_ex2_results} shows the posterior histograms and the
trace plots of the Markov chain samples generated by the algorithm. 

It is to be noted that similar to the first numerical example, we
observed wrong convergence of the samples using the standard Gaussian likelihood function.

\subsection{A two-dimensional material inversion problem}

This example concerns the inference of wave speed in a two dimensional
spatial domain, consisting of ten non-overlapping sub-domains, each
with a different unknown wave speed.

\begin{figure}[!h]
\begin{center}
\includegraphics[width=0.4\linewidth]{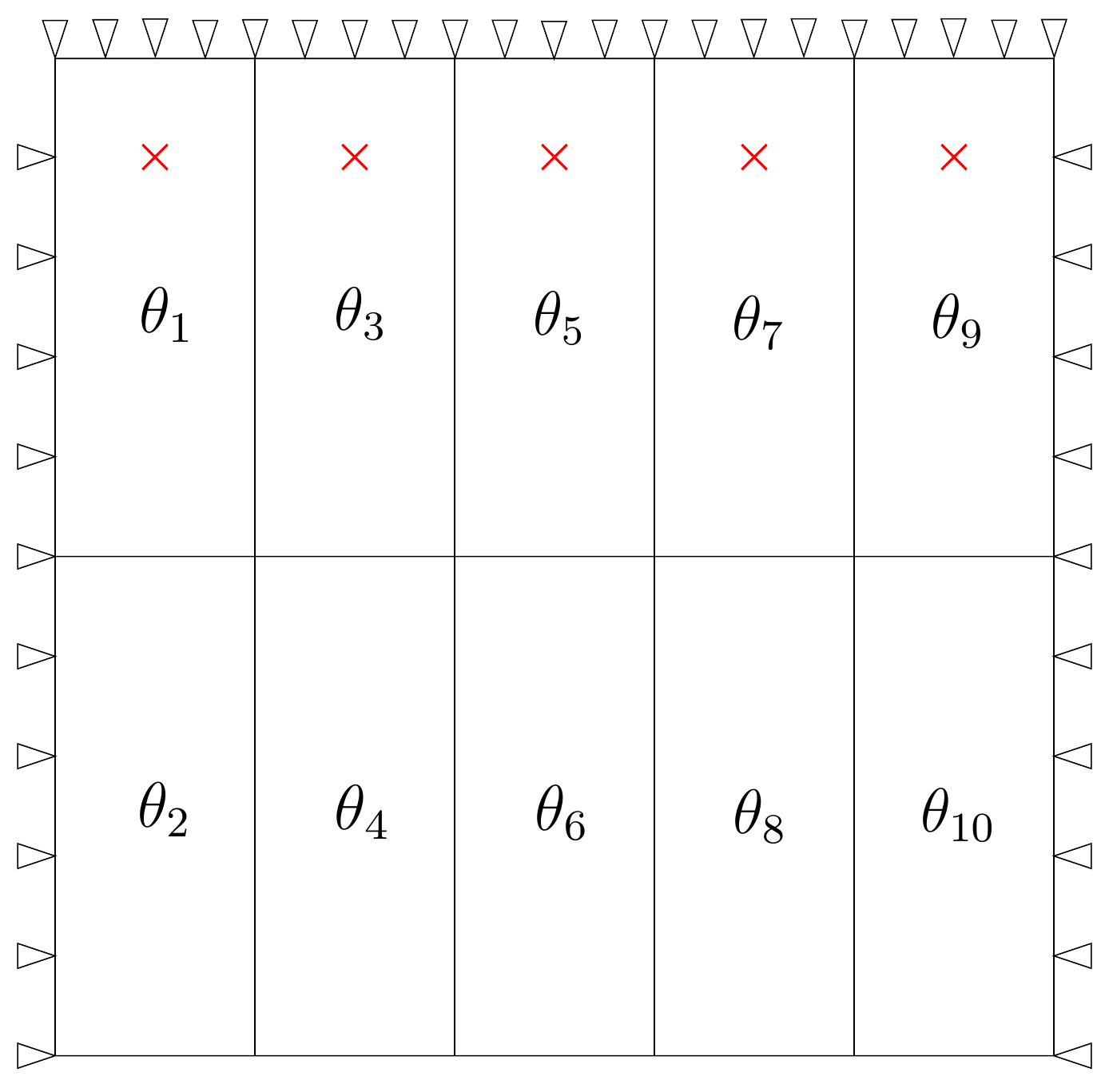}     
\vspace{-.01cm}
\caption{A schematic representation of the computational domain, the
  wave speed structure, and the arrays of sources (red $\times$) and
  the uniform array of receivers (black $\triangledown$).}
\label{EX3_domain_setting}
\end{center}
\end{figure}

\medskip
\noindent
{\bf Problem formulation.} 
Consider the initial-boundary value problem for the acoustic wave equation
in a two-dimensional spatial domain $D=[-1 ,1]\times[-2,0]$:
\begin{equation}\label{wave_2D}
\begin{array}{ll}
u_{tt}(t,{\bf x}) - \nabla_{\bf x} \cdot ( a^2({\bf x};\boldsymbol\theta) \,
  \nabla_{\bf x} u(t,{\bf x}) = F_s(t, {\bf x}), & \ \ \ t \in [0,T],
  \quad {\bf x} \in D,\\
u(t,{\bf x})=0, \qquad u_t(t,{\bf x})=0, & \ \ \ t=0,
  \quad {\bf x} \in D\\
u(t,{\bf x}) = 0, & \ \ \ t \in [0,T],
  \quad {\bf x} \in \partial D.
\end{array}    
\end{equation}
The source $F_s$ is a Ricker wavelet applied on a square
$D_s$ of size $0.08 \times 0.08$ centered at a given point
$(x_{1,s}, x_{2,s})$:
$$
F_s(t,{\bf x}) = \left\{  \begin{array}{ll}
10^3 \, (1-20 \, (t-0.1)^2) \, \exp(-10 \, (t-0.1)^2), & \ \ \ {\bf x} \in D_s,\\
0, & \ \ \ {\bf x} \not\in D_s.
\end{array}\right.
$$
The wave speed $a$ is assumed to be described by ten (unknown)
parameters $\boldsymbol\theta=(\theta_1,\dotsc,\theta_{10})$. The
computational domain $D$ is split into ten non-overlapping
sub-domains of size $0.4 \times 1$, each with a constant velocity
$a=\theta_i$, with $i=1, \dotsc, 10$; see Figure \ref{EX3_domain_setting}.

\medskip
\noindent
{\bf Synthetic data.} The acquisition geometry consists of $N_s=5$
sources $F_s$ evenly centered at depth $x_2 = - 0.2$ along the horizontal
line $x_1 \in [-0.8, 0.8]$, and $N_r=201$ evenly spaced receivers at
the surface ($x_2 = 0$) with a spacing of 0.02 and the left ($x_1 =
-1$) and right ($x_1 =1$) boundaries with a spacing of 0.04. 
Figure \ref{EX3_domain_setting} shows the computational domain and the
configuration of the problem. 

Each receiver located at ${\bf x}_r$, with $r=1, \dotsc, N_r$, records a
(noisy) discrete-time signal $g_s(t_k,{\bf x}_r)$, due to a single
source $F_s$, with $s=1, \dotsc, N_s$, and over the time interval
$[0,T]$ at $N$ discrete
time levels $t_k = (k-1) \, \Delta t$, with $\Delta t = T/(N - 1)$ and
$k=1, \dotsc, N$. Moreover, let $f_s(t_k,{\bf x}_r; \boldsymbol\theta)$ denote the
corresponding simulated signal for a given $\boldsymbol\theta$,
obtained by employing a finite difference scheme based on
the second-order central difference discretization of \eqref{wave_2D}; see
\cite{Motamed_etal:13} for details of the deterministic solver. 
We use a uniform grid with spatial grid-lengths $\Delta x_1 =
\Delta x_2 = 0.02$ and time step $\Delta t = 0.005$. 
We set the final time $T=4$ with $N=801$ time steps. 
We choose a fixed parameter vector 
$$
\boldsymbol\theta^*=(3, 2, 3.5, 2.5,4,3,4.5,3.5,5,4),
$$
and generate synthetic data by first computing $f_s(t_k,{\bf x}_r;
\boldsymbol\theta^*)$ and then polluting it with a noise of the form \eqref{noise_structure}:
$$
g_s(t_k, {\bf x}_r) = \varepsilon_{skr}^{(1)}  \, f_s(t_k, {\bf x}_r; \boldsymbol\theta^*) +
\varepsilon_{skr}^{(2)}, \quad \varepsilon_{skr}^{(1)} \sim
\text{Gamma}(1000, 1000), \quad \varepsilon_{skr}^{(2)} \sim
\text{Unif}(-0.05, 0.05).
$$
Let ${\bf g}_{sr} \in {\mathbb R}^{N}$ and ${\bf
  f}_{sr}(\boldsymbol\theta) \in {\mathbb R}^{N}$ denote the vectors of
recorded and simulated discrete-time signals due to a single source
$F_s$, $s=1, \dotsc, N_s$, at the receiver location
${\bf x}_r$, $r=1, \dotsc, N_r$. Note that since there are $N_s =5$ single sources, we have in total $N_s N_r=1005$ such signals.

\medskip
\noindent
{\bf Computations.} The goal is to employ Bayesian inversion and compute the
conditional posterior of the model parameter vector $\boldsymbol\theta
= (\theta_1, \dotsc,\theta_{10})$. 
We follow Algorithm \ref{ALG_MG} with the
following choices:

\begin{itemize}
\item {\it Likelihood}: 
$
\pi_{\text{exp}}( {\bf g} | \boldsymbol\theta) =s^N \, \exp ( -s
\, \sum_{s=1}^{N_s}\, \sum_{r=1}^{N_r} d_W({\bf f}_{sr}(\boldsymbol\theta), {\bf g}_{sr})).
$

\item {\it Priors}: 
$\theta_j \sim \text{Unif}(1 ,6),  \ \ j=1, \dotsc, 10, \qquad 
s \sim \text{Gamma}(1 ,0.1).
$

\item {\it Proposal}: A Gaussian random walk
  \eqref{proposal_random_walk} with covariance
  $\Sigma$, which is initially chosen to be diagonal and
  then updated based on the first 1000 MCMC samples.

\item {\it Initial point}: $\theta_j^{(0)} = 3.8, \ \ j=1, \dotsc, 10,
  \qquad s^{(0)}=5000$.

\end{itemize}

\begin{figure}[!h]
\center
\includegraphics[width=0.99\linewidth]{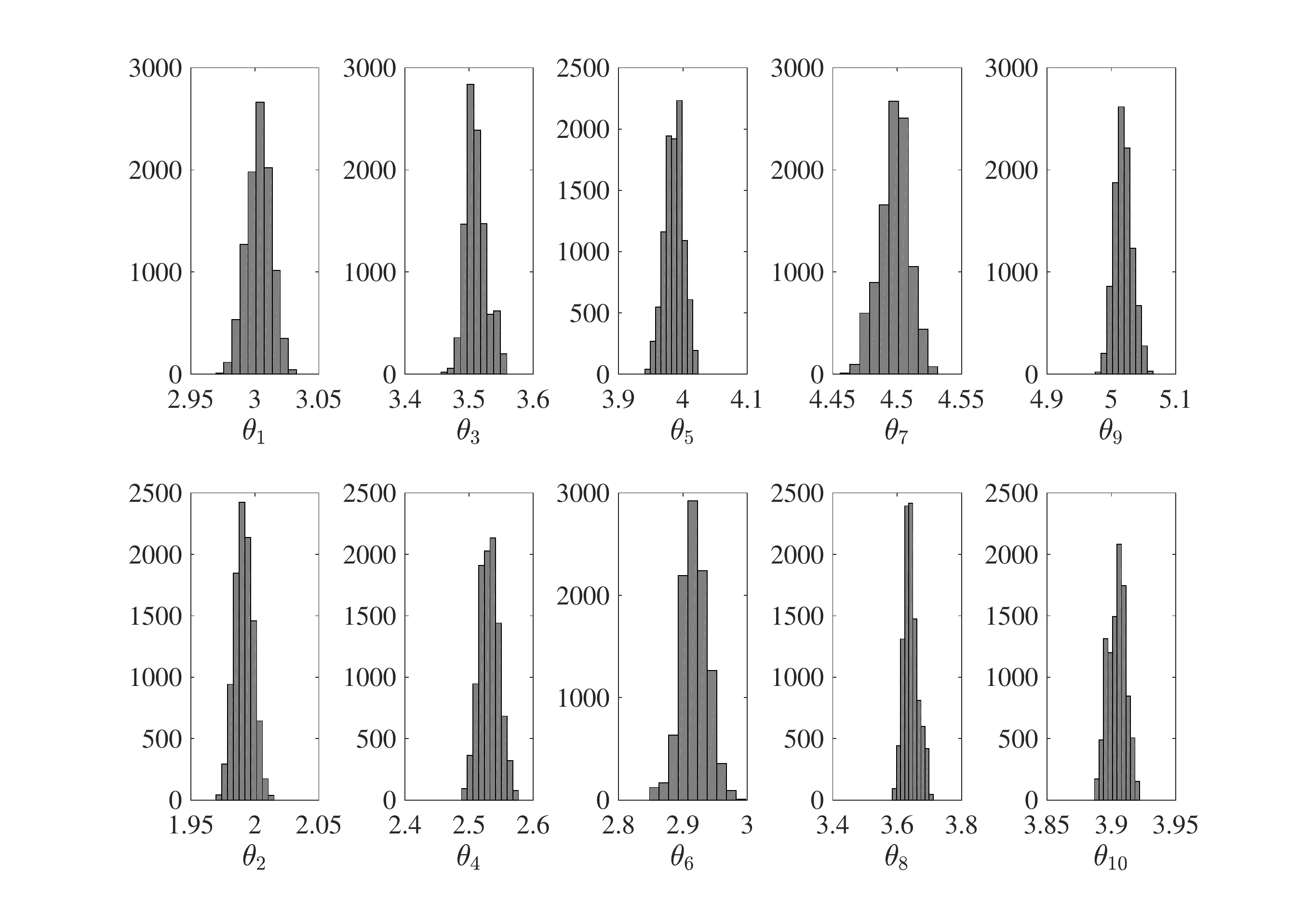}
\vspace{-.1cm}
\caption{Posterior histograms of the Markov chain
  samples generated by the proposed Wasserstein-driven Bayesian algorithm.} 
\label{fig_ex3_results}
\end{figure}

We run the algorithm with $M_0=80000$ iterations and remove the first
$M_b=65000$ samples. We also use a
thinning period of $M_t=3$. This would give a total of $M=5000$ Markov
chain samples. 
Figure \ref{fig_ex3_results} shows the posterior
histograms of the Markov chain samples generated by the proposed algorithm. 

\begin{figure}[!h]
\center
\subfigure[$\boldsymbol\theta^*$]{\includegraphics[width=0.3\linewidth]{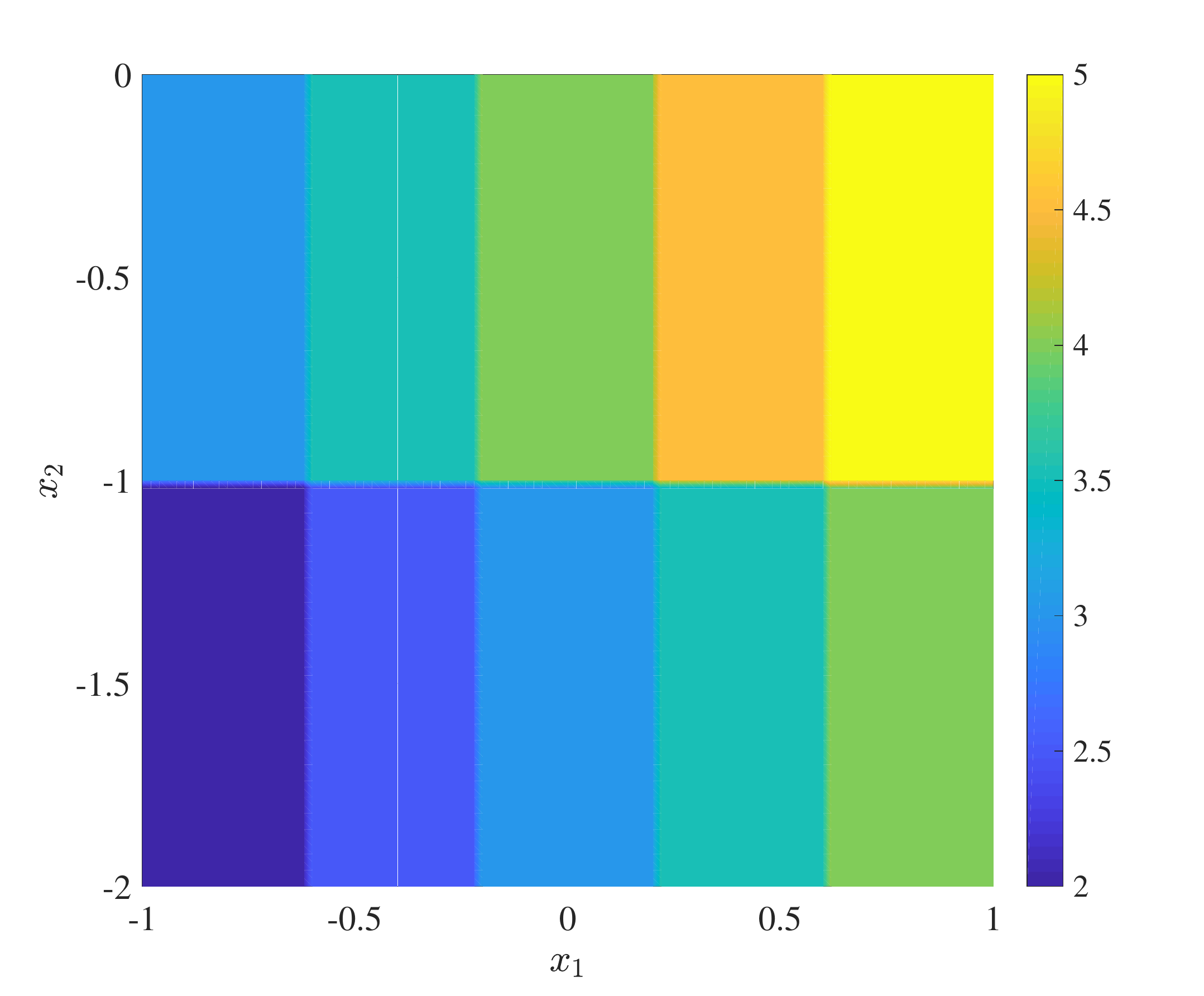}}
\hskip .01cm
\subfigure[{${\mathbb E}[\boldsymbol\theta]$}]{\includegraphics[width=0.3\linewidth]{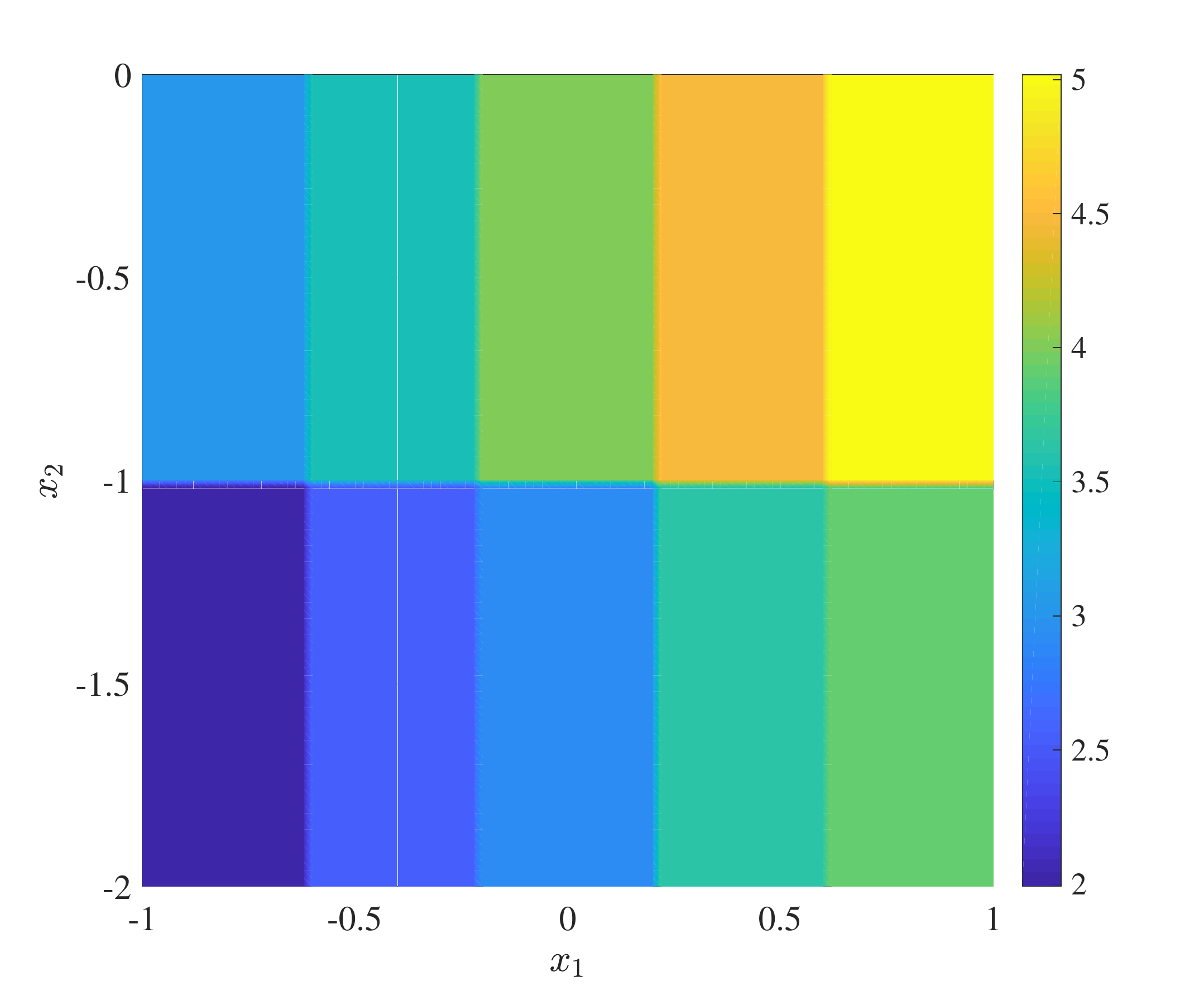}}
\hskip .01cm
\subfigure[{${\mathbb V}[\boldsymbol\theta]^{1/2}$}]{\includegraphics[width=0.3\linewidth]{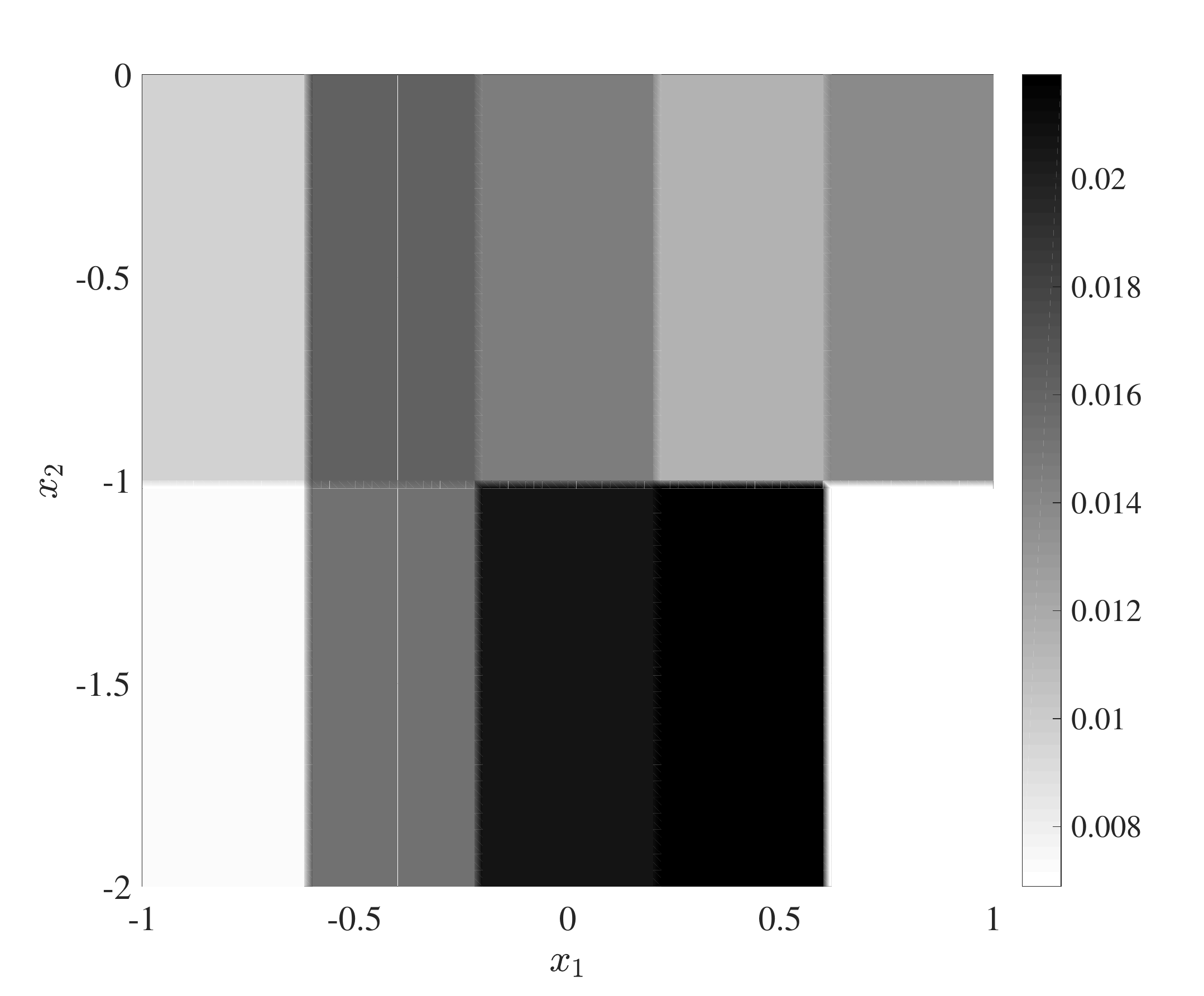}}
\caption{Map of the original noise-free parameters (a), the mean map
  (b) and the standard deviation map (c) obtained by the generated
  MCMC samples of the parameter posteriors.} 
\label{fig_ex3_vel}
\end{figure}

Figure \ref{fig_ex3_vel} shows a map of the original noise-free
parameters (a), the mean map (b) and the standard deviation map (c)
obtained by the generated MCMC samples of the parameter posteriors.

\vspace{-.2cm}

\section{Conclusion}

We have presented a robust Bayesian framework based on a new exponential
likelihood function driven by the quadratic Wasserstein metric. We
demonstrated that this framework is able to treat complicated noise
structures and presented several numerical examples subject to combined additive and multiplicative noise. As in the deterministic setting we observe that the convexity properties of the Wasserstein metric translates into better convergence, here towards the  correct posteriors. 



{\color{black}It is to be noted that in this work we exclusively used a
  linear  ``shift and rescale`` approach to convert signals into
  positive probabilities. Although this procedure was shown to produce
  good results, there may be other alternatives, e.g. the combination
  of the exponential and linear procedure as discussed in the end of
  Section 3, that could prove to be more suitable for certain
  problems. A detailed exploration of the impact of normalization on the
  proposed method is a subject of our future work.}

As the framework presented here is particularly suitable for wave
propagation problems we also plan to apply it to geoacoustic inversion
of water and seafloor parameters in shallow water range dependent
environments. There the material properties vary smoothly and
consequently the forward problems can be efficiently simulated with
our arbitrary order wave solvers \cite{secondHermite}. 

{\color{black}Other future directions include the application of Wasserstein
distance to other Bayesian models, such as variational inference (see
e.g. \cite{VI:17}), that
suggest alternative strategies to MCMC sampling and hence tend to be
faster to scale to large data.}

\bigskip
{\bf Acknowledgements.} 
The authors are indebted to Dr. Gabriel Huerta at Sandia National
Laboratories for his guidance and stimulating discussions. We also thank the reviewers for their constructive comments.   

\bibliographystyle{plain}
\bibliography{refs}

\end{document}